\documentclass[reqno, draft]{amsart}
\usepackage{amsmath,amsthm,amscd,amssymb,latexsym,upref}

\newcommand{\bbN}{{\mathbb{N}}}
\newcommand{\bbR}{{\mathbb{R}}}
\newcommand{\bbP}{{\mathbb{P}}}
\newcommand{\bbZ}{{\mathbb{Z}}}
\newcommand{\bbC}{{\mathbb{C}}}

\newcommand{\bbQ}{{\mathbb{Q}}}

\newcommand{\calB}{{\mathcal B}}
\newcommand{\cB}{{\mathcal B}}
\newcommand{\calC}{{\mathcal C}}
\newcommand{\calD}{{\mathcal D}}
\newcommand{\calE}{{\mathcal E}}
\newcommand{\calF}{{\mathcal F}}

\newcommand{\cH}{{\mathcal H}}
\newcommand{\calK}{{\mathcal K}}
\newcommand{\calL}{{\mathcal L}}
\newcommand{\calM}{{\mathcal M}}

\newcommand{\calP}{{\mathcal P}}

\newcommand{\calV}{{\mathcal V}}

\renewcommand{\gg}{\mathfrak{g}}

\newcommand{\hatt}{\widehat}  
\newcommand{\Pinf}{P_\infty}
\newcommand{\Div}{\operatorname{Div}}

\newcommand{\no}{\nonumber}
\newcommand{\lb}{\label}
\newcommand{\f}{\frac}
\newcommand{\ul}{\underline}
\newcommand{\ol}{\overline}
\newcommand{\ti}{\tilde}
\newcommand{\wti}{\widetilde}
\newcommand{\N}{n}
\newcommand{\uc}{{\underline{c}}}
\newcommand{\pa}{\partial}
\newcommand{\Oh}{O}
\newcommand{\oh}{o}

\newcommand{\loc}{\text{\rm{loc}}}

\newcommand{\ran}{\text{\rm{ran}}}

\newcommand{\dom}{\text{\rm{dom}}}
\newcommand{\bi}{\bibitem}
\newcommand{\humu}{{ \hat{\underline{\mu} }}}

\newcommand{\hmu}{{\hat{\mu} }}

\newcommand{\uz}{{\underline{z}}}

\newcommand{\uxi}{{\underline{\Xi}}}
\newcommand{\ome}{\omega}

\newcommand{\al}{\alpha}
\newcommand{\ual}{{\underline{\alpha}}}
\newcommand{\ua}{{\underline{A}}}
\newcommand{\hua}{{ \underline{\hatt{A} }}}

\newcommand{\uU}{{\underline{U}}}

\newcommand{\Mod}{\text{\rm{Mod}}}
\newcommand{\Sp}{\text{\rm Sp}}

\renewcommand{\Re}{\text{\rm Re}}
\renewcommand{\Im}{\text{\rm Im}}

\DeclareMathOperator{\diag}{diag}

\DeclareMathOperator{\tr}{tr}
\DeclareMathOperator{\KdV}{KdV}
\DeclareMathOperator{\sKdV}{s-KdV}
\newcommand\shKdV{{\mathop{\text{\rm s-}\widehat{\text{\rm KdV}}}}}

\DeclareMathOperator{\sym}{Sym}
\newcommand{\symn}{{\sym^n (\calK_n)}}

\allowdisplaybreaks
\numberwithin{equation}{section}

\newtheorem{theorem}{Theorem}[section]
\newtheorem{lemma}[theorem]{Lemma}

\theoremstyle{definition}
\newtheorem{definition}[theorem]{Definition}
\newtheorem{hypothesis}[theorem]{Hypothesis}
\newtheorem{remark}[theorem]{Remark}
\newtheorem{example}[theorem]{Example}

\begin{document}
\title[The spectrum of quasi-periodic algebro-geometric 
KdV potentials]{On the spectrum of Schr\"odinger operators
with quasi-periodic algebro-geometric  KdV potentials}
\author[V.\ Batchenko]{Volodymyr Batchenko}
\address{Department of Mathematics,
University of Missouri,
Columbia, MO 65211, USA}
\email{batchenv@math.missouri.edu}
\author[F.\ Gesztesy]{Fritz Gesztesy}
\address{Department of Mathematics,
University of Missouri,
Columbia, MO 65211, USA}
\email{fritz@math.missouri.edu}
\urladdr{http://www.math.missouri.edu/people/fgesztesy.html}
\dedicatory{Dedicated with great pleasure to Vladimir A.\ Marchenko on the
occasion of his 80th birthday.}
\date{\today}
\begin{abstract}
We characterize the spectrum of one-dimensional Schr\"odinger
operators $H=-d^2/dx^2+V$ in $L^2(\bbR;dx)$ with quasi-periodic
complex-valued algebro-geometric potentials $V$ (i.e., potentials $V$
which satisfy one (and hence infinitely many) equation(s) of the
stationary Korteweg--de Vries (KdV) hierarchy) associated with nonsingular
hyperelliptic curves. The corresponding problem appears to have been open
since the mid-seventies. The spectrum of $H$ coincides with the
conditional stability set of $H$ and can explicitly be described in terms
of the mean value of the inverse of the diagonal Green's function of $H$.

As a result, the spectrum of $H$ consists of finitely many simple
analytic arcs and one semi-infinite simple analytic arc in the
complex plane. Crossings as well as confluences of spectral arcs are
possible and discussed as well. These results extend to the 
$L^p(\bbR;dx)$-setting for $p\in [1,\infty)$. 
\end{abstract}

\maketitle

\section{Introduction}\lb{s1}

It is well-known since the work of Novikov \cite{No74}, Marchenko
\cite{Ma74}, \cite{Ma74a}, Dubrovin \cite{Du75}, Dubrovin, Matveev, and
Novikov \cite{DMN76}, Flaschka \cite{Fl75}, Its and Matveev \cite{IM75}, 
Lax \cite{La75}, McKean and van Moerbeke \cite{MM75} (see also
\cite[Sects.\ 3.4, 3.5]{BBEIM94}, \cite[p.\ 111--112, App.\ J]{GH03},
\cite[Sect.\ 4.4]{Ma86}, \cite[Sects.\ II.6--II.10]{NMPZ84} and the
references therein) that the self-adjoint Schr\"odinger operator 
\begin{equation}
H=-\f{d^2}{dx^2}+V, \quad \dom(H)=H^{2,2}(\bbR) \lb{1.1}
\end{equation}
in $L^2(\bbR;dx)$ with a real-valued periodic, or more generally, {\it
quasi-periodic} and {\it real-valued} potential $V$, that satisfies one
(and hence infinitely many) equation(s) of the stationary Korteweg--de
Vries (KdV) equations, leads to a finite-gap, or perhaps more
appropriately, to a finite-band spectrum $\sigma (H)$ of the form
\begin{equation}
\sigma(H)=\bigcup_{m=0}^{n-1} [E_{2m},E_{2m+1}]\cup [E_{2n},\infty).
\lb{1.2}
\end{equation}
It is also well-known, due to work of Serov \cite{Se60} and
Rofe-Beketov \cite{Ro63} in 1960 and 1963, respectively (see also
\cite{Tk64}), that if $V$ is {\it periodic} and {\it complex-valued} then
the spectrum of the non-self-adjoint Schr\"odinger operator $H$ defined as
in \eqref{1.1} consists either of infinitely many simple analytic arcs,
or else, of a finite number of simple analytic arcs and one
semi-infinite simple analytic arc tending to infinity. It seems
plausible that the latter case is again connected with (complex-valued)
stationary solutions of equations of the KdV hierarchy, but to the best
of our knowledge, this has not been studied in the literature. In
particular, the next scenario in line, the determination of the
spectrum of $H$ in the case of {\it quasi-periodic} and {\it
complex-valued} solutions of the stationary KdV equation apparently has
never been clarified. The latter problem is open since the mid-seventies
and it is the purpose of this paper to provide a comprehensive solution of
it.

To describe our results, a bit of preparation is needed. Let 
\begin{equation}
G(z,x,x')=(H-z)^{-1}(x,x'), \quad z\in\bbC\backslash\sigma(H), \; 
x,x'\in\bbR,  \lb{1.1a}
\end{equation} 
be the Green's function of $H$ (here $\sigma(H)$ denotes the spectrum
of $H$) and denote by $g(z,x)$ the corresponding diagonal
Green's function of $H$ defined by
\begin{align}
&g(z,x)=G(z,x,x) =\f{i\prod_{j=1}^n [z-\mu_j(x)]}{2R_{2n+1}(z)^{1/2}}, 
\lb{1.2a} \\
& R_{2n+1}(z)=\prod_{m=0}^{2n} (z-E_m), \quad 
\{E_m\}_{m=0}^{2n}\subset\bbC,  \lb{1.2b} \\
& E_m\neq E_{m'} \text{  for $m\neq m'$, \; 
$m,m'=0,1,\dots,2n$.} \lb{1.2c}
\end{align}
For any quasi-periodic (in fact, Bohr (uniformly) almost periodic)
function $f$ the mean value $\langle f\rangle$ of $f$ is defined by
\begin{equation}
\langle f\rangle =\lim_{R\to\infty}\f{1}{2R} \int_{-R}^{R} dx \,
f(x). \lb{1.3}
\end{equation}
Moreover, we introduce the set $\Sigma$ by
\begin{equation}
\Sigma=\big\{\lambda\in\bbC\,\big|\, \Re\big(\big\langle
g(\lambda,\cdot)^{-1}\big\rangle\big)=0\big\} \lb{1.4}
\end{equation}
and note that
\begin{equation}
\langle g(z,\cdot)\rangle=\f{i\prod_{j=1}^n
\big(z-\wti\lambda_j\big)}{2R_{2n+1}(z)^{1/2}} \lb{1.4a}
\end{equation}
for some constants $\{\wti\lambda_j\}_{j=1}^n\subset\bbC$.

Finally, we denote by $\sigma_{\rm p}(T)$, 
$\sigma_{\rm r}(T)$, $\sigma_{\rm c}(T)$, 
$\sigma_{\rm{e}}(T)$, and $\sigma_{\rm{ap}}(T)$, the point spectrum
(i.e., the set of eigenvalues), the residual spectrum, the continuous
spectrum, the essential spectrum (cf.\ \eqref{5.22b}), and the
approximate point spectrum of a densely defined closed operator
$T$ in a complex Hilbert space, respectively.

Our principal new results, to be proved in Section \ref{s4}, then read
as follows:

\begin{theorem}  \lb{t1.1} 
Assume that $V$ is a quasi-periodic $($complex-valued\,$)$ solution of the
$n$th stationary KdV equation associated with the hyperelliptic curve
$y^2=R_{2n+1}(z)$ subject to \eqref{1.2b} and \eqref{1.2c}. Then the
following assertions hold:
\\
$(i)$ The point spectrum and residual spectrum of $H$ are empty and
hence the spectrum of $H$ is purely continuous,
\begin{align}
&\sigma_{\rm p}(H)=\sigma_{\rm r}(H)=\emptyset,  \lb{1.5} \\
&\sigma(H)=\sigma_{\rm c}(H)=\sigma_{\rm e}(H)=\sigma_{\rm ap}(H).
\lb{1.6}
\end{align}
$(ii)$ The spectrum of $H$ coincides with $\Sigma$ and equals the
conditional stability set of $H$,
\begin{align}
\sigma(H) &=\big\{\lambda\in\bbC\,\big|\, \Re\big(\big\langle
g(\lambda,\cdot)^{-1}\big\rangle\big)=0\big\}  \lb{1.10} \\ 
&=\{\lambda\in\bbC\,|\, \text{there exists at least one bounded
distributional solution}  \no \\
& \hspace*{1.8cm} \text{$0\neq\psi\in L^\infty(\bbR;dx)$ of
$H\psi=\lambda\psi$}\}.  \lb{1.11} 
\end{align}
$(iii)$ $\sigma(H)$ is contained in the semi-strip
\begin{equation}
\sigma(H)\subset \{z\in\bbC\,|\, \Im(z)\in [M_1,M_2], \, \Re(z)\geq
M_3\},  \lb{1.12}
\end{equation}
where
\begin{equation}
M_1=\inf_{x\in\bbR}[\Im(V(x))], \quad 
M_2=\sup_{x\in\bbR}[\Im(V(x))], \quad M_3=
\inf_{x\in\bbR}[\Re(V(x))]. \lb{1.13}
\end{equation}
$(iv)$ $\sigma(H)$ consists of finitely many simple analytic arcs and
one simple semi-infinite arc. These analytic arcs may only end at the
points $\wti\lambda_1,\dots,\wti\lambda_n$, $E_0,\dots,E_{2n}$, and
at infinity. The semi-infinite arc, $\sigma_\infty$, asymptotically
approaches the half-line 
$L_{\langle V\rangle}=\{z\in\bbC \,|\, z=\langle V\rangle +x, \, x\geq
0\}$ in the following sense: asymptotically, 
$\sigma_\infty$ can be parameterized by
\begin{equation}
\sigma_\infty=\big\{z\in\bbC \,\big|\, z=R+i\,\Im(\langle
V\rangle) +\Oh\big(R^{-1/2}\big) 
\text{ as $R\uparrow\infty$}\big\}.  \lb{1.16}
\end{equation} 
$(v)$ Each $E_m$, $m=0,\dots,2n$, is met by 
at least one of these arcs. More precisely, a particular $E_{m_0}$ is
hit by precisely $2N_0+1$ analytic arcs, where $N_0\in\{0,\dots,n\}$
denotes the number of $\wti\lambda_j$ that coincide with $E_{m_0}$.
Adjacent arcs meet at an angle $2\pi/(2N_0+1)$ at $E_{m_0}$. $($Thus,
generically, $N_0=0$ and precisely one arc hits $E_{m_0}$.$)$ \\
$(vi)$ Crossings of spectral arcs are permitted and take place
precisely when 
\begin{equation}
\Re\big(\big\langle g(\wti\lambda_{j_0},\cdot)^{-1}\big\rangle\big)=0 
\, \text{ for some $j_0\in\{1,\dots,n\}$ with $\wti\lambda_{j_0}\notin
\{E_m\}_{m=0}^{2n}$}. \lb{1.17}
\end{equation}
In this case $2M_0+2$ analytic arcs are converging toward
$\wti\lambda_{j_0}$, where $M_0\in\{1,\dots,n\}$ denotes the number of
$\wti\lambda_j$ that coincide with $\wti\lambda_{j_0}$. Adjacent arcs meet
at an angle $\pi/(M_0+1)$ at $\wti\lambda_{j_0}$. $($Thus, generically,
$M_0=1$ and two arcs cross at a right angle.$)$ \\
$(vii)$ The resolvent set $\bbC\backslash\sigma(H)$ of $H$ is
path-connected.
\end{theorem}

Naturally, Theorem \ref{t1.1} applies to the special case where $V$ is 
a periodic (complex-valued) solution of the $n$th stationary KdV
equation associated with a nonsingular hyperelliptic curve. Even in this
special case, items (v) and (vi) of Theorem \ref{t1.1} provide additional
new details on the nature of the spectrum of $H$. 

As described in Remark \ref{r4.10}, these results extend to the
$L^p(\bbR;dx)$-setting for $p\in [1,\infty)$.

Theorem \ref{t1.1} focuses on stationary quasi-periodic solutions of the
KdV hierarchy for the following reasons. First of all, the class of
algebro-geometric solutions of the (time-dependent) KdV hierarchy is
defined as the class of all solutions of some (and hence infinitely
many) equations of the stationary KdV hierarchy. Secondly,
time-dependent algebro-geometric solutions of a particular equation of
the (time-dependent) KdV hierarchy just represent isospectral
deformations (the deformation parameter being the time variable) of a
fixed stationary algebro-geometric KdV solution (the latter can be
viewed as the initial condition at a fixed time $t_0$). In the present  
case of quasi-periodic algebro-geometric solutions of the $n$th KdV
equation, the isospectral manifold of such a given solution is a complex 
$n$-dimensional torus, and time-dependent solutions trace out a path
in that isospectral torus (cf.\ the discussion in \cite[p.\ 12]{GH03}).

Finally, we give a brief discussion of the contents of each section.
In Section \ref{s2} we provide the necessary background material
including a quick construction of the KdV hierarchy of nonlinear
evolution equations and its Lax pairs using a polynomial recursion
formalism. We also discuss the hyperelliptic Riemann surface
underlying the stationary KdV hierarchy, the corresponding
Baker--Akhiezer function, and the necessary ingredients to describe
the Its--Matveev formula for stationary KdV solutions. Section
\ref{s3} focuses on the diagonal Green's function of the Schr\"odinger
operator $H$, a key ingredient in our characterization of the spectrum
$\sigma(H)$ of $H$ in Section \ref{s4} (cf.\ \eqref{1.10}). Our
principal Section \ref{s4} is then devoted to a proof of Theorem
\ref{t1.1}. Appendix \ref{sA} provides the necessary summary of tools
needed from elementary algebraic geometry (most notably the theory of
compact (hyperelliptic) Riemann surfaces) and sets the stage for some
of the notation used in Sections \ref{s2}--\ref{s4}. Appendix \ref{sB}
provides additional insight into one ingredient of the Its--Matveev
formula; Appendix \ref{sC} illustrates our results in the special
periodic non-self-adjoint case and provides a simple yet nontrivial
example in the elliptic genus one case. 

Our methods extend to the case of algebro-geometric
non-self-adjoint second order finite difference (Jacobi) operators
associated with the Toda lattice hierarchy. Moreover, they extend to
the infinite genus limit $n\to\infty$ (cf.\ \eqref{1.2}--\eqref{1.2b})
using the approach in \cite{Ge01}. This will be studied elsewhere.

\smallskip

{\bf Dedication.} 
It is with great pleasure that we dedicate this paper to Vladimir A.\
Marchenko on the occasion of his 80th birthday. His strong influence on
the subject at hand is universally admired.   

\section{The KdV hierarchy, hyperelliptic curves, \\ and the Its--Matveev
formula} \label{s2}

In this section we briefly review the recursive construction of the KdV 
hierarchy and associated Lax pairs following \cite{GRT96} and especially, 
\cite[Ch.\ 1]{GH03}. Moreover, we discuss the class of algebro-geometric
solutions of the KdV hierarchy corresponding to the underlying
hyperelliptic curve and recall the Its--Matveev formula for such
solutions. The material in this preparatory section is known and detailed
accounts with proofs can be found, for instance, in \cite[Ch.\ 1]{GH03}.
For the notation employed in connection with elementary concepts in
algebraic geometry (more precisely, the theory of compact Riemann
surfaces), we refer to Appendix \ref{sA}. 

Throughout this section we suppose the hypothesis
\begin{equation}
V\in C^\infty (\bbR)  \lb{1.2.1a}
\end{equation}
and consider the one-dimensional Schr\"odinger differential expression
\begin{equation}
L = - \frac{d^2}{dx^2} + V. \lb{1.2.2} 
\end{equation}
To construct the KdV hierarchy we need a second differential expression
$P_{2n+1}$ of order $2n+1$, $n\in\bbN_0$, defined recursively in the
following. We take the quickest  route to the construction of $P_{2n+1}$,
and hence to that of the KdV hierarchy, by starting from the recursion
relation \eqref{1.2.3} below. 

Define $\{f_\ell\}_{\ell\in\bbN_0}$ recursively by
\begin{equation}
f_0 = 1,\quad f_{\ell,x} = - (1/4) f_{\ell-1, xxx} + Vf_{\ell-1,x}
+ (1/2) V_x f_{\ell-1},\quad \ell\in\bbN.
\lb{1.2.3}
\end{equation}
Explicitly, one finds
\begin{align}
f_0 & =1, \no \\
f_1 & = \tfrac12 V + c_1,  \no \\
f_2 &= - \tfrac18 V_{xx} +\tfrac38 V^2 +
 c_1\tfrac12 V + c_2, \lb{1.2.4} \\ 
f_3 & = \tfrac1{32} V_{xxxx} -
\tfrac5{16} VV_{xx} - \tfrac5{32} V_x^2 +
\tfrac5{16} V^3\no \\
&\qquad + c_1 \big(- \tfrac18 V_{xx}
+ \tfrac38 V^2 \big)+c_2 \tfrac12 V + c_3, \quad\text{etc.} \no
\end{align}
Here $\{c_k\}_{k\in\bbN}\subset\bbC$ denote integration 
constants which naturally arise when solving \eqref{1.2.3}. 

Subsequently, it will be convenient to also introduce 
the corresponding homogeneous coefficients $\hat f_\ell$, defined by
the vanishing of the integration constants $c_k$ for $k=1,\dots,\ell$, 
\begin{equation}
\hat{f}_0=f_0=1, \quad \hat{f}_\ell=f_\ell\big|_{c_k=0, \,
k=1,\dots,\ell}, \quad \ell\in\bbN. \lb{1.2.4a}
\end{equation}
Hence,
\begin{equation}
f_\ell=\sum_{k=0}^{\ell}c_{\ell-k}\hat{f}_{k}, \quad \ell\in\bbN_0, \lb{1.2.4c}
\end{equation}
introducing
\begin{equation}
c_0=1. \lb{1.2.4ca}
\end{equation}
One can prove inductively that all homogeneous elements $\hat f_\ell$ (and
hence all $f_\ell$) are differential polynomials in $V$, that is,
polynomials with respect to $V$ and its $x$-derivatives up to order
$2\ell-2$, $\ell\in\bbN$. 

Next we define differential expressions
$P_{2n+1}$  of order $2n+1$ by 
\begin{equation}
P_{2n+1}  = \sum_{\ell=0}^n \Big( f_{n-\ell}
\frac{d}{dx}-\frac12 f_{n-\ell,x} \Big) L^{\ell}, \quad n\in\bbN_0.
\lb{1.2.5}
\end{equation}
Using the recursion \eqref{1.2.3}, the commutator of $P_{2n+1}$ and $L$
can be explicitly computed and one obtains 
\begin{equation}
[P_{2n+1},L] = 2f_{n+1,x}, \quad n\in\bbN_0.
\lb{1.2.6}
\end{equation}
In particular, $(L,P_{2n+1})$ 
represents the celebrated \textit{Lax pair} of the KdV hierarchy. Varying
$n\in\bbN_0$, the stationary KdV  hierarchy is then defined in terms of
the vanishing of the commutator of $P_{2n+1}$ and $L$ in \eqref{1.2.6}
by\footnote{In a slight abuse of notation we 
will occasionally stress the functional dependence of $f_\ell$ 
on $V$, writing $f_\ell(V)$.},
\begin{equation}
-[P_{2n+1}, L] =-2f_{n+1,x}(V)=\sKdV_{n}(V)=0,\quad n\in\bbN_0.
\lb{1.2.7}
\end{equation}
Explicitly, 
\begin{align} 
\sKdV_{0}(V)&=-V_x =0, \no \\
\sKdV_{1}(V)&=\tfrac14 V_{xxx} -\tfrac32 VV_x + c_1(-V_x)=0,\lb{1.2.8}\\
\sKdV_{2}(V)&=-\tfrac1{16} V_{xxxxx} +\tfrac58 V_{xxx} + \tfrac54
V_x V_{xx} - \tfrac{15}8 V^2 V_x \no \\ 
&\quad +c_1 \big( \tfrac14 V_{xxx} -\tfrac32 VV_x\big)+ c_2(-V_x) =0,
\quad \text{etc.,} \no
\end{align}
represent the first few equations of the stationary KdV hierarchy. By 
definition, the set of solutions of \eqref{1.2.7}, with $n$ ranging 
in $\bbN_0$ and $c_k$ in $\bbC$, $k\in\bbN$, represents the 
class of algebro-geometric KdV solutions. At times it will be convenient
to  abbreviate algebro-geometric stationary KdV solutions $V$ simply as
KdV \textit{potentials}. 

In the following we will frequently assume that $V$ satisfies the $n$th
stationary KdV equation. By this we mean it satisfies one of the
$n$th stationary KdV equations after a particular choice of integration
constants $c_k\in\bbC$, $k=1,\dots,n$, $n\in\bbN$, has been made. 

Next, we introduce a polynomial $F_n$ of degree $n$ with respect to 
the spectral parameter $z\in\bbC$ by
\begin{equation}
F_n (z,x)  = \sum_{\ell=0}^n f_{n-\ell}(x) z^\ell.  \lb{1.2.9}
\end{equation}
Explicitly, one obtains
\begin{align}
F_0&=1, \no \\
F_1&=z+\tfrac12 V + c_1, \no \\
F_2&=z^2+\tfrac12 Vz- \tfrac18 V_{xx} +\tfrac38 V^2 +
 c_1\big(\tfrac12 V+z\big) + c_2, \label{1.2.9a} \\
F_3&=z^3+\tfrac12 Vz^2+\big(- \tfrac18 V_{xx} +
\tfrac38 V^2\big)z  
+\tfrac1{32} V_{xxxx} -\tfrac5{16} VV_{xx} -
\tfrac5{32} V_x^2  \no \\
& \quad+\tfrac5{16} V^3+c_1\big(z^2+\tfrac12 Vz-
\tfrac18 V_{xx} +\tfrac38 V^2\big)+c_2\big(z+\tfrac12 V\big)+c_3,
\quad \text{etc.} \no
\end{align}
The recursion relation \eqref{1.2.3} and equation \eqref{1.2.7} imply that
\begin{equation}
F_{n,xxx}  -4 (V-z) F_{n,x}  -2V_x F_n  =0. \lb{1.2.10}
\end{equation}
Multiplying \eqref{1.2.10} by $F_n$, a subsequent integration with respect
to $x$ results in 
\begin{equation}
(1/2) F_{n,xx}  F_n  - (1/4) F_{n,x}^2  -(V-z) F_n^2  = R_{2n+1},
\lb{1.2.11}
\end{equation}
where $R_{2n+1}$ is a monic polynomial of degree $2n+1$.  We denote its
roots by $\{E_m\}_{m=0}^{2n}$, and hence  write
\begin{equation}
R_{2n+1} (z) = \prod_{m=0}^{2n} (z-E_m),
\quad \{ E_m\}_{m=0}^{2n} \subset \bbC.
\lb{1.2.12}
\end{equation}
One can show that equation \eqref{1.2.11} leads to an explicit
determination of the integration constants $c_1,\dots,c_n$ in 
\begin{equation}
\sKdV_n(V)=-2f_{n+1,x}(V)=0 \lb{1.2.12A}
\end{equation}
in terms of the zeros $E_0,\dots,E_{2n}$ of the associated
polynomial $R_{2n+1}$ in \eqref{1.2.12}. In fact, one can prove 
\begin{equation}
c_k=c_k(\ul E), \quad k=1,\dots,n, \lb{1.2.12B}
\end{equation}
where 
\begin{align}
c_k(\ul E)&=-\!\!\!\!\!\sum_{\substack{j_0,\dots,j_{2n}=0\\
 j_0+\cdots+j_{2n}=k}}^{k}\!\!
\f{(2j_0)!\cdots(2j_{2n})!}
{2^{2k} (j_0!)^2\cdots (j_{2n}!)^2 (2j_0-1)\cdots(2j_{2n}-1)}
E_0^{j_0}\cdots E_{2n}^{j_{2n}}, \no \\
& \hspace*{7.5cm} k=1,\dots,n. \label{1.2.12C} 
\end{align}

\begin{remark} \lb{r2.1}
Suppose $V\in C^{2n+1}(\bbR)$ satisfies the $n$th stationary KdV
equation $\sKdV_{n}(V)=-2f_{n+1,x}(V)=0$ for a given set of integration
constants $c_k$, $k=1,\dots,n$. Introducing
$F_n$ as in \eqref{1.2.9} with $f_0,\dots,f_n$ given by \eqref{1.2.4c}
then yields equation \eqref{1.2.10} and hence \eqref{1.2.11}. The latter
equation in turn, as shown inductively in \cite[Prop.\ 2.1]{GW96}, yields  
\begin{equation}
V\in C^\infty(\bbR) \, \text{ and } \, f_\ell \in C^\infty(\bbR), \;
\ell=0,\dots,n. \lb{2.20}
\end{equation} 
Thus, without loss of generality, we may assume in the following that
solutions of $\sKdV_{n}(V)=0$ satisfy $V\in C^\infty(\bbR)$.
\end{remark}

Next, we study the restriction of the differential expression $P_{2n+1}$
to the two-dimensional  kernel (i.e., the formal null space in an
algebraic sense as opposed to the functional analytic one) of $(L-z)$.
More precisely, let
\begin{equation}
\ker (L -z)=\left\{ \psi\colon \bbR\to\bbC_\infty \text{ meromorphic} 
\mid (L-z)\psi=0 \right\}, \quad z\in\bbC. \lb{1.2.13}
\end{equation}
Then \eqref{1.2.5} implies
\begin{equation}
P_{2n+1}\big|_{\ker(L-z)} =
\Big(F_n (z) \frac{d}{dx} -\frac12F_{n,x}(z)\Big)\Big|_{\ker(L-z)}.
\lb{1.2.14}
\end{equation}
We emphasize that the result \eqref{1.2.14} is valid independently of whether or
not $P_{2n+1}$ and $L$ commute.  However, if one makes the additional
assumption that $P_{2n+1}$ and $L$ commute, one can prove that this implies
an algebraic relationship between $P_{2n+1}$ and $L$. 
\begin{theorem} \label{t1.2.2}
Fix $n\in\bbN_0$ and assume that $P_{2n+1}$ and $L$ commute,
$[P_{2n+1},L]=0$, or equivalently, suppose 
$\sKdV_{n}(V)=-2f_{n+1,x}(V)=0$. Then $L$ and $P_{2n+1}$ satisfy an
algebraic relationship of the type $($cf.\ \eqref{1.2.12}$)$
\begin{align}
\begin{split}
&\calF_n(L,-iP_{2n+1})= - P_{2n+1}^2-R_{2n+1}(L)=0, \label{1.2.15} \\
& R_{2n+1}(z) = \prod_{m=0}^{2n} (z-E_m), \quad z\in \bbC. 
\end{split} 
\end{align}
\end{theorem}

The expression $\calF_n(L,-iP_{2n+1})$ is called the Burchnall--Chaundy
polynomial of the  pair $(L,P_{2n+1})$. Equation \eqref{1.2.15} naturally
leads to  the hyperelliptic curve $\calK_n$ of (arithmetic) genus
$n\in\bbN_0$ (possibly with a singular affine part), where
\begin{align}
\begin{split}
&\calK_n \colon \calF_n(z,y)=y^2-R_{2n+1}(z)=0,  \lb{1.2.17} \\
&R_{2n+1}(z) = \prod_{m=0}^{2n} (z-E_m), \quad 
\{ E_m\}_{m=0}^{2n} \subset \bbC. 
\end{split}
\end{align}
The curve $\calK_n$ is compactified by joining the point $P_\infty$ but for
notational simplicity  the compactification is also denoted by $\calK_n$.
Points $P$ on $\calK_n\backslash\{P_\infty\}$ are  represented as pairs
$P=(z,y)$, where $y(\cdot)$ is the meromorphic function on $\calK_n$
satisfying $\calF_n(z,y)=0$. The complex structure on 
$\calK_n$ is then defined in the usual way, see Appendix \ref{sA}. 
Hence, $\calK_n$ becomes a two-sheeted hyperelliptic Riemann
surface of (arithmetic) genus $n\in\bbN_0$ (possibly with a singular
affine part) in a standard manner.

We also emphasize that by fixing the curve $\calK_n$ (i.e., by fixing
$E_0,\dots,E_{2n}$), the integration constants $c_1,\dots,c_n$ in
$f_{n+1,x}$ (and hence in the corresponding stationary $\KdV_n$ equation)
are uniquely determined as is clear from \eqref{1.2.12B} and 
\eqref{1.2.12C}, which establish the integration constants $c_k$ as
symmetric functions of $E_0,\dots,E_{2n}$. 

For notational simplicity we will usually tacitly assume that $n\in\bbN$. 
The trivial case $n=0$ which leads to $V(x)=E_0$ is of no interest to us
in this paper. 

In the following, the zeros\footnote{If $V\in L^\infty(\bbR;dx)$, these
zeros (generically) are the Dirichlet  eigenvalues of a closed operator in
$L^2(\bbR)$ associated with the differential expression $L$ and a
Dirichlet boundary condition  at $x\in\bbR$.} of the polynomial
$F_n(\cdot,x)$ (cf.\ \eqref{1.2.9}) will play a special role. We denote
them by $\{\mu_j(x)\}_{j=1}^n$ and hence write 
\begin{equation}
F_n(z,x) =\prod_{j=1}^n [z-\mu_j(x)].
\lb{1.3.2}
\end{equation}
{}From \eqref{1.2.11} we see that
\begin{equation}
R_{2n+1} + (1/4) F_{n,x}^2=F_{n}H_{n+1},
\lb{1.3.3}
\end{equation}
where
\begin{equation}
H_{n+1}(z,x)=(1/2) F_{n,xx}(z,x)+(z-V(x))F_n(z,x) \lb{1.3.4}
\end{equation}
is a monic polynomial of degree $n+1$. We introduce the corresponding  
roots\footnote{If $V\in L^\infty(\bbR;dx)$, these roots (generically) are
the Neumann  eigenvalues of a  closed operator in $L^2(\bbR)$ associated
with $L$ and a Neumann boundary  condition at $x\in\bbR$.}
$\{\nu_\ell(x)\}_{\ell=0}^n$ of $H_{n+1}(\cdot,x)$ by
\begin{equation}
H_{n+1}(z,x) =\prod_{\ell=0}^n [z-\nu_\ell(x)]. \lb{1.3.6}
\end{equation}
Explicitly, one computes from \eqref{1.2.4} and \eqref{1.2.9},
\begin{align}
H_1 &  = z-V, \no \\
H_2 & = z^2-\tfrac12 Vz+ \tfrac14 V_{xx} -
\tfrac12 V^2 + c_1(z-V), \lb{1.3.5} \\
H_3 & = z^3- \tfrac12 V z^2+\tfrac18\big(V_{xx} -V^2\big)z
-\tfrac1{16} V_{xxxx} +\tfrac38 V_x^2
+ \tfrac12 VV_{xx} \no \\
& \quad -\tfrac38 V^3  +c_1\big(z^2-\tfrac12 Vz+ \tfrac14
V_{xx} -\tfrac12 V^2\big)+c_2(z-V),\quad  \text{etc.} \no
\end{align}
The next step is crucial; it permits us to ``lift'' the zeros $\mu_j$ and
$\nu_\ell$ of $F_n$ and $H_{n+1}$ from $\bbC$ to the curve $\calK_n$. 
{}From \eqref{1.3.3} one infers 
\begin{equation} 
R_{2n+1}(z) + (1/4) F_{n,x}(z)^2=0, \quad 
z\in\{\mu_j,\nu_\ell\}_{j=1,\dots,n, \ell=0,\dots,n}. \lb{1.3.7}
\end{equation}
We now introduce
$\{ \hat \mu_j(x) \}_{j=1,\dots,n}\subset \calK_n$ and 
$\{ \hat \nu_\ell (x) \}_{\ell=0,\dots,n}\subset \calK_n$ by
\begin{equation}
\hat \mu_j(x)=(\mu_j(x),-(i/2) F_{n,x}(\mu_j(x), x)), \quad j=1,\dots,n, \,
x\in\bbR \lb{1.3.8a}
\end{equation}
and
\begin{equation}
\hat \nu_\ell(x)=(\nu_\ell(x),(i/2) F_{n,x}(\nu_\ell(x), x)),
\quad \ell=0,\dots,n, \, x\in\bbR. \lb{1.3.8b}
\end{equation}
Due to the $C^\infty (\bbR)$ assumption \eqref{1.2.1a} on $V$, 
$F_n(z,\cdot)\in C^\infty(\bbR)$ by \eqref{1.2.3} and \eqref{1.2.9}, and 
hence also $H_{n+1}(z,\cdot)\in C^\infty(\bbR)$ by \eqref{1.3.4}. Thus,
one concludes 
\begin{equation}
\mu_j, \nu_\ell \in C(\bbR), \;\, j=1,\dots,n, \; \ell=0,\dots,n, 
\lb{1.3.8c}
\end{equation}
taking multiplicities (and appropriate renumbering) of the zeros of $F_n$
and $H_{n+1}$ into account. (Away from collisions of zeros, $\mu_j$ and
$\nu_\ell$ are of course $C^\infty$.)

Next, we define the fundamental meromorphic function
$\phi(\cdot,x)$ on $\calK_n$,
\begin{align}
\phi (P,x) &= \frac{iy + (1/2) F_{n,x} (z, x)}{F_n (z, x)} \lb{1.3.9a}\\
& = \frac{-H_{n+1} (z,x)}{iy - (1/2) F_{n,x}(z, x)}, \lb{1.3.9b} \\ 
&\, P = (z,y)\in\calK_n, \; x\in\bbR \no
\end{align}
with divisor $(\phi(\cdot, x))$ of $\phi(\cdot, x)$ given by
\begin{equation}
(\phi (\cdot, x))=\calD_{\hat \nu_0(x) \hat{\underline \nu} (x)} -
\calD_{P_\infty \hat{\underline \mu} (x)}, \lb{1.3.10}
\end{equation}
using \eqref{1.3.2}, \eqref{1.3.6}, and \eqref{1.3.8c}. Here we
abbreviated
\begin{equation}
\hat{\underline \mu} = \{\hat \mu_1, \dots, \hat \mu_n\}, \, 
\hat{\underline \nu} = \{\hat \nu_1, \dots,\hat \nu_n\}\in\symn 
\lb{1.3.11}
\end{equation}
(cf.\ the notation introduced in Appendix \ref{sA}). 
The stationary Baker--Akhiezer function $\psi(\cdot, x, x_0)$ on
$\calK_n\backslash \{ P_\infty \}$ is then defined in terms of
$\phi(\cdot,x)$ by
\begin{equation}
\psi(P,x,x_0) =\exp \bigg( \int_{x_0}^x  dx' \, \phi(P,x')\bigg), \quad
P\in\calK_n\backslash \{ P_\infty \},\, (x,x_0) \in\bbR^2.
\lb{1.3.12}
\end{equation}
Basic properties of $\phi$ and $\psi$ are summarized
in the following result (where $W(f,g)=fg^\prime-f^\prime g$ denotes the
Wronskian of $f$ and $g$, and $P^*$ abbreviates $P^*=(z,-y)$ for
$P=(z,y)$).
\begin{lemma} \lb{l1.3.1}  
Assume $V\in C^\infty(\bbR)$ satisfies the $n$th stationary KdV
equation \eqref{1.2.7}. Moreover, let 
$P= (z,y) \in \calK_n \backslash \{P_\infty\}$ and   
$(x,x_0) \in \bbR^2$.  Then $\phi$ satisfies
the Riccati-type equation
\begin{align}
& \phi_x (P) + \phi(P)^2 = V -z, \lb{1.3.14}
\intertext{as well as}
& \phi(P) \phi (P^*)= \f{H_{n+1} (z)}{F_n (z)},\lb{1.3.16}\\
&\phi (P) + \phi (P^*)= \f{F_{n,x} (z)}{F_n (z)},\lb{1.3.17}\\
& \phi(P)-\phi (P^*)=\f{2iy}{F_n (z)}. \lb{1.3.18}
\end{align}
Moreover, $\psi$ satisfies 
\begin{align}
&(L-z(P)) \psi(P) =0, \quad (P_{2n+1}-iy(P)) \psi(P) =0, \lb{1.3.15} \\
&\psi (P,x,x_0) =\bigg(\f{F_n (z,x)}{F_n (z,x_0)} \bigg)^{1/2} \exp
\bigg(iy \int_{x_0}^x  dx' \, F_n (z,x')^{-1} \bigg), \lb{1.3.20} \\
&\psi (P,x,x_0) \psi (P^*, x, x_0)= \f{F_n (z,x)}{F_n(z,x_0)},\lb{1.3.19}\\
&\psi_x (P,x,x_0) \psi_x (P^*, x, x_0)= \f{H_{n+1} (z,x)}{F_n (z,x_0)},
\lb{1.3.19a}\\
&\psi(P,x,x_0)\psi_x(P^*,x,x_0)
+\psi(P^*,x,x_0)\psi_x(P,x,x_0)=\f{F_{n,x}(z,x)}{F_n(z,x_0)}, 
\lb{1.3.19aaa} \\
&W(\psi (P,\cdot,x_0), \psi (P^*, \cdot, x_0))=-\f{2iy}{F_{n}(z,x_0)}. 
\lb{1.3.19aab} 
\end{align}
In addition, as long as the zeros of $F_n(\cdot,x)$ are all simple for 
$x\in\Omega$, $\Omega\subseteq\bbR$ an open interval, $\psi(\cdot,x,x_0)$
is meromorphic on $\calK_n\backslash\{P_\infty\}$ for $x,x_0\in\Omega$.
\end{lemma} 

Combining the polynomial recursion approach with \eqref{1.3.2} 
readily yields trace formulas for the KdV invariants, that is, expressions 
of $f_\ell$ in terms of symmetric functions of the zeros $\mu_j$ of $F_n$. 

\begin{lemma} \lb{l1.3.9a}
Assume $V\in C^\infty(\bbR)$ satisfies the $n$th stationary KdV
equation \eqref{1.2.7}. Then, 
\begin{align}
V&=\sum_{m=0}^{2n} E_m -2\sum_{j=1}^n \mu_j, \lb{1.3.59a} \\
V^2-(1/2) V_{xx}&=\sum_{m=0}^{2n} E_m^2 -2\sum_{j=1}^n \mu_j^2, 
\text{ etc.}\lb{1.3.59b}
\end{align}
\end{lemma}

Equation \eqref{1.3.59a} represents the trace formula for the
algebro-geometric potential $V$. In addition, \eqref{1.3.59b}
indicates that higher-order trace formulas associated with the KdV
hierarchy can be obtained from \eqref{1.3.2} comparing powers of
$z$. We omit further details and refer to \cite[Ch.\ 1]{GH03} and
\cite{GRT96}.

From this point on we assume that the affine part of $\calK_n$ is
nonsingular, that is, 
\begin{equation}
E_m\neq E_{m'} \text{  for $m\neq m'$, \; $m,m'=0,1,\dots,2n$}.\lb{1.3.52A}
\end{equation}

Since nonspecial divisors play a fundamental role in this context we also
recall the following fact.

\begin{lemma} \lb{l1.3.9ba}
Suppose that the affine part of $\calK_n$ is nonsingular and 
assume that $V\in C^\infty(\bbR)\cap L^\infty(\bbR;dx)$ satisfies the
$n$th stationary KdV equation \eqref{1.2.7}. Let $\calD_{\humu}$,
$\humu=(\hmu_1,\dots,\hmu_n)$ be the Dirichlet divisor of degree
$n$ associated with $V$ defined according to \eqref{1.3.8a}, that is,
\begin{equation}
\hmu_j(x)=(\mu_j(x),-(i/2)F_{n,x}(\mu_j(x),x)), \quad j=1,\dots,n, \;
x\in\bbR. \lb{1.3.59AA}
\end{equation}
Then $\calD_{\humu(x)}$ is nonspecial for all $x\in\bbR$. Moreover, there
exists a constant $C>0$ such that
\begin{equation}
|\mu_j(x)|\leq C, \quad j=1,\dots,n, \; x\in\bbR. \lb{1.3C}
\end{equation}
\end{lemma}

\begin{remark} \lb{r2.7} 
Assume that $V\in C^\infty(\bbR)\cap L^\infty(\bbR;dx)$ satisfies the
$n$th stationary KdV equation \eqref{1.2.7}. We recall that 
$f_\ell\in C^\infty(\bbR)$, $\ell\in\bbN_0$, by \eqref{2.20} since
$f_\ell$ are differential polynomials in $V$. Moreover, we note that
\eqref{1.3C} implies that  $f_\ell\in L^\infty(\bbR;dx)$,
$\ell=0,\dots,n$, employing the fact that $f_\ell$, $\ell=0,\dots,n$, are
elementary symmetric functions of $\mu_1,\dots,\mu_n$ (cf.\ \eqref{1.2.9}
and \eqref{1.3.2}). Since $f_{n+1,x}=0$, one can use the recursion
relation \eqref{1.2.3} to reduce $f_k$ for  $k\geq n+2$ to a linear
combination of $f_1,\dots,f_n$. Thus, 
\begin{equation} 
f_\ell\in C^\infty(\bbR)\cap L^\infty(\bbR;dx), \quad \ell\in\bbN_0.
\lb{2.57}
\end{equation}
Using the fact that for fixed $1\leq p \leq \infty$,
\begin{equation}
h, h^{(k)} \in L^p(\bbR;dx) \, \text{ imply } \, h^{(\ell)} \in
L^p(\bbR;dx), \;\; \ell=1,\dots,k-1 \lb{2.58}
\end{equation}
(cf., e.g., \cite[p.\ 168--170]{BB83}), one then infers
\begin{equation}
V^{(\ell)} \in L^\infty(\bbR;dx), \quad \ell\in\bbN_0, \lb{2.59}
\end{equation}
applying \eqref{2.58} with $p=\infty$.
\end{remark}

We continue with the theta function representation for $\psi$ and $V$.
For general background information and the notation employed we refer to
Appendix \ref{sA}. 

Let $\theta$ denote the Riemann theta function associated with
$\calK_n$ (whose affine part is assumed to be nonsingular) and a fixed
homology basis $\{a_j,b_j\}_{j=1}^n$ on $\calK_n$. Next, choosing a base
point $Q_0\in\calK_n\backslash\Pinf$, the Abel maps
$\ua_{Q_0}$ and $\ual_{Q_0}$ are  defined by \eqref{aa46} and \eqref{aa47},
and the Riemann vector $\underline{\Xi}_{Q_0}$ is given by \eqref{aa55}.

Next, let $\ome_{P_\infty,0}^{(2)}$ denote the normalized differential of 
the second kind defined by
\begin{align}
&\ome_{P_\infty,0}^{(2)}  = -\f1{2 y} \prod_{j=1}^n(z -\lambda_j) d z
\underset{\zeta\to 0}{=}\big(\zeta^{-2}+\Oh(1)\big)d\zeta \text{  as $P\to
P_\infty$},  \lb{1.3.56j1} \\
& \hspace*{5.6cm} \zeta=\sigma /z^{1/2}, \; \sigma\in\{1,-1\}, \no 
\end{align}
where the constants $\lambda_j\in\bbC$, $j=1,\dots, n$, are determined by
employing the normalization
\begin{equation}
\int_{a_j}\ome_{P_\infty,0}^{(2)}=0, \quad j=1,\dots, n. \lb{1.3.71a}
\end{equation}
One then infers
\begin{equation}
\int_{Q_0}^P \ome_{P_\infty,0}^{(2)} \underset{\zeta\to 0}{=}-\zeta^{-1}
+ e^{(2)}_0(Q_0)+\Oh(\zeta) 
\text{ as $P\to P_\infty$} \lb{1.3.71b} 
\end{equation}
for some constant $e^{(2)}_0(Q_0)\in\bbC$. The vector of $b$-periods of
$\ome_{P_\infty,0}^{(2)}/(2\pi i)$ is denoted by  
\begin{equation}
\uU_0^{(2)} = (U_{0,1}^{(2)},
\dots, U_{0,n}^{(2)} ),\quad U_{0,j}^{(2)}
=\f{1}{2\pi i} \int_{b_j} \ome_{P_\infty,0}^{(2)}, \; j=1,\dots,n.
\lb{1.3.60b}
\end{equation}
By \eqref{b27b} one concludes
\begin{equation}
U_{0,j}^{(2)}=-2 c_j(n), \quad j=1,\dots,n. \lb{1.3.60c}
\end{equation}
In the following it will be convenient to introduce the abbreviation
\begin{equation}
\uz(P,\ul Q) =\ul\Xi_{Q_0}-\ul A_{Q_0}(P)
+\ul\alpha_{Q_0}(\calD_{\ul Q}),  \;
P\in\calK_\N, \; \ul Q=\{Q_1,\dots, Q_\N\} \in \sym^\N (\calK_\N). 
\lb{1.3.60d} 
\end{equation}
We note that $\ul{z}(\cdot,\ul{Q})$ is independent of the choice of base
point $Q_0$.

\begin{theorem} \lb{t1.3.10} 
Suppose that $V\in C^\infty(\bbR)\cap L^\infty(\bbR;dx)$ satisfies the
$n$th stationary KdV equation \eqref{1.2.7} on $\bbR$. In addition,
assume the affine part of $\calK_n$ to be nonsingular and let $P\in
\calK_n \backslash \{ P_\infty\}$ and $x,x_0\in\bbR$. Then 
$\calD_{\ul{\hat\mu}(x)}$ and $\calD_{\ul{\hat\nu} (x)}$ are 
nonspecial for $x\in\bbR$. Moreover,\footnote{To avoid multi-valued 
expressions in formulas 
such as \eqref{1.3.59}, etc., we agree to always 
choose the same path of integration  connecting $Q_0$ and $P$ and 
refer to Remark \ref{raa26a} for additional tacitly  assumed
conventions.} 
\begin{align}
\psi(P,x,x_0) &= \frac{\theta (\uz(P_\infty,\humu (x_0)))\theta
(\uz(P,\humu (x)))}{\theta (\uz(P_\infty,\humu (x)))\theta (\uz(P,\humu
(x_0)))} \no \\ 
& \quad \times \exp \bigg[-i (x-x_0)\bigg(\int_{Q_0}^P
\ome_{P_\infty,0}^{(2)}-e^{(2)}_0(Q_0)\bigg)\bigg], \lb{1.3.59}
\end{align}
with the linearizing property of the Abel map,
\begin{equation}
 \ual_{Q_0} (\calD_{\humu(x)}) = \Big(
\ual_{Q_0} (\calD_{\humu(x_0)}) + i\uU_0^{(2)}(x-x_0)\Big) \pmod {L_n}.
\lb{1.3.60} 
\end{equation}
The Its--Matveev formula for $V$ reads
\begin{equation}
 V(x) = E_0 + \sum_{j=1}^n (E_{2j-1} +E_{2j} - 2\lambda_j) 
-2\pa_x^2 \ln\big(\theta (\uxi_{Q_0} -\ua_{Q_0} (P_\infty) +
\ual_{Q_0} (\calD_{\humu (x)}))\big). \lb{1.3.61a}
\end{equation}
\end{theorem}

Combining \eqref{1.3.60} and \eqref{1.3.61a} shows the remarkable linearity
of the theta function with respect to $x$ in the Its--Matveev formula for
$V$. In fact, one can rewrite \eqref{1.3.61a} as 
\begin{equation}
V(x)=\Lambda_0 -2\partial_x^2 \ln(\theta(\ul A +\ul B x)),  \lb{1.3.IM} 
\end{equation}
where 
\begin{align}
\ul A&= \ul \Xi_{Q_0}-\ua_{Q_0} (P_\infty)-i\uU_0^{(2)}x_0 
+ \ual_{Q_0}(\calD_{\humu (x_0)}),  \lb{1.3.IMA} \\
\ul B&=i\uU_0^{(2)}, \lb{1.3.IMB} \\
\Lambda_0&=E_0+ \sum_{j=1}^n (E_{2j-1}+ E_{2j} - 2\lambda_j). \lb{1.3.IML}
\end{align}
Hence the constants $\Lambda_0\in\bbC$ and $\ul B \in\bbC^n$ are
uniquely determined by $\calK_n$ (and its homology basis), and the
constant $\ul A\in\bbC^n$ is in one-to-one correspondence with the
Dirichlet data $\humu(x_0)=(\hmu_1(x_0),\dots,\hmu_n(x_0))
\in\symn$ at the point $x_0$.

\begin{remark} \lb{r2.8}
If one assumes $V$ in \eqref{1.3.61a} $($or \eqref{1.3.IM}$)$ to be
quasi-periodic (cf.\ \eqref{3.14a} and \eqref{3.14b}), then there exists a
homology basis $\{\tilde a_j, \tilde b_j\}_{j=1}^n$ on $\calK_n$ such
that 
$\wti{\ul B}=i\wti{\ul U}^{(2)}_0$ satisfies the constraint
\begin{equation}
\wti{\ul B}=i \wti{\ul U}^{(2)}_0 \in \bbR^n. \lb{2.70}
\end{equation}
This is studied in detail in Appendix \ref{sB}.
\end{remark}

An example illustrating some of the general results of this section is
provided in Appendix \ref{sC}.

\section{The diagonal Green's function of $H$}
\lb{s3}

In this section we focus on the diagonal Green's function of $H$ and
derive a variety of results to be used in our principal Section
\ref{s4}.  

We start with some preparations. We denote by 
\begin{equation}
W(f,g)(x)=f(x)g_x(x)-f_x(x)g(x) \, \text{ for a.e. $x\in\bbR$} \lb{3.1}
\end{equation}
the Wronskian of $f,g \in AC_{\loc}(\bbR)$ (with $AC_{\loc}(\bbR)$ the set
of locally absolutely continuous functions on $\bbR$).
\begin{lemma} \lb{l3.1}
Assume\footnote{One could admit more severe local singularities; in
particular, one could assume $q$ to be meromorphic, but we will not need
this in this paper.} $q\in L^1_{\loc}(\bbR)$, define $\tau=-d^2/dx^2+q$,
and let $u_j(z)$, $j=1,2$ be two $($not necessarily distinct$)$
distributional  solutions\footnote{That is, $u, u_x \in AC_{\loc}(\bbR)$.}
of $\tau u=zu$ for some $z\in\bbC$. Define $U(z,x)=u_1(z,x)u_2(z,x)$,
$(z,x)\in\bbC\times\bbR$. Then,
\begin{equation}
2U_{xx}U-U_x^2-4(q-z)U^2=-W(u_1,u_2)^2. \lb{3.3}
\end{equation} 
If in addition $q_x\in L^1_{\loc}(\bbR)$, then 
\begin{equation}
U_{xxx}-4(q-z)U_x-2q_xU=0. \lb{3.2} 
\end{equation} 
\end{lemma}
\begin{proof}
Equation \eqref{3.2} is a well-known fact going back to at least  
Appell \cite{Ap80}. Equation \eqref{3.3} either follows upon integration
using the integrating factor $U$, or alternatively, can be verified
directly from the definition of $U$. We omit the straightforward
computations. 
\end{proof}

Introducing
\begin{equation}
\gg(z,x)=u_1(z,x)u_2(z,x)/W(u_1(z),u_2(z)), \quad z\in\bbC, \; x\in\bbR, 
\lb{3.4}
\end{equation}
Lemma \ref{l3.1} implies the following result. 

\begin{lemma} \lb{l3.2}
Assume that $q\in L^1_{\loc}(\bbR)$ and $(z,x)\in\bbC\times\bbR$. Then,
\begin{align}
&2\gg_{xx}\gg-\gg_x^2-4(q-z)\gg^2=-1, \lb{3.6} \\
&-\big(\gg^{-1}\big)_z=2\gg  
 +\big\{\gg\big[u_1^{-2}W(u_1,u_{1,z}) +
u_2^{-2}W(u_2,u_{2,z})\big]\big\}_x, \lb{3.7} \\
&-\big(\gg^{-1}\big)_z=2\gg - \gg_{xxz}+\big[\gg^{-1}\gg_x\gg_z\big]_x
\lb{3.8} \\
& \hspace*{1.48cm} =2\gg - 
\Big\{\Big[\big(\gg^{-1}\big)\big(\gg^{-1}\big)_{zx}
-\big(\gg^{-1}\big)_x \big(\gg^{-1}\big)_z
\Big]\Big/\big(\gg^{-3}\big)\Big\}_x. \lb{3.8a} 
\end{align} 
If in addition $q_x\in L^1_{\loc}(\bbR)$, then 
\begin{equation}
\gg_{xxx}-4(q-z)\gg_x-2q_x\gg=0. \lb{3.5} 
\end{equation}
\end{lemma}
\begin{proof}
Equations \eqref{3.5} and \eqref{3.6} are clear from \eqref{3.2} and
\eqref{3.3}. Equation \eqref{3.7} follows from
\begin{equation}
\big(\gg^{-1})_z=u_2^{-2}W(u_2,u_{2,z})-u_1^{-2}W(u_1,u_{1,z}) \lb{3.9}
\end{equation}
and 
\begin{equation}
W(u_j,u_{j,z})_x=-u_j^2, \quad j=1,2. \lb{3.10}
\end{equation}
Finally, \eqref{3.8a} (and hence \eqref{3.8}) follows from \eqref{3.4},
\eqref{3.6}, and
\eqref{3.7} by a straightforward, though tedious, computation.  
\end{proof}
\noindent Equation \eqref{3.8} is known and can be found, for instance, in
\cite{GD75}. Similarly, \eqref{3.7} can be inferred, for example, from the
results in \cite[p.\ 369]{CL90}. 

Next, we turn to the analog of $\gg$ in connection with the
algebro-geometric potential $V$ in \eqref{1.3.61a}. Introducing 
\begin{equation}
g(P,x)=\f{\psi(P,x,x_0)\psi(P^*,x,x_0)}{W(\psi(P,\cdot,x_0),
\psi(P^*,\cdot,x_0))}, \quad P\in \calK_n
\backslash \{ P_\infty\}, \; x,x_0\in\bbR, \lb{3.11}
\end{equation}
equations \eqref{1.3.19} and \eqref{1.3.19aab} imply
\begin{equation}
g(P,x)=\f{i F_n(z,x)}{2y}, \quad P=(z,y)\in \calK_n
\backslash \{ P_\infty\}, \; x\in\bbR. \lb{3.12}
\end{equation}
Together with $g(P,x)$ we also introduce its two branches $g_\pm(z,x)$
defined on the upper and lower sheets $\Pi_\pm$ of $\calK_n$ (cf.\
\eqref{b2}, \eqref{b3}, and \eqref{b22})
\begin{equation}
g_\pm (z,x)=\pm \f{i F_n(z,x)}{2R_{2n+1}(z)^{1/2}}, \quad z\in\Pi, \;
x\in\bbR \lb{3.13}
\end{equation}
with $\Pi=\bbC\backslash\calC$ the cut plane introduced in \eqref{b3}. A 
comparison of \eqref{3.4}, \eqref{3.11}--\eqref{3.13}, then shows that
$g_{\pm}(z,\cdot)$ satisfy \eqref{3.6}--\eqref{3.5}. 

For convenience we will subsequently focus on $g_+$ whenever possible and
then use the simplified notation
\begin{equation}
g(z,x)=g_+(z,x), \quad z\in\Pi, \; x\in\bbR. \lb{3.14}
\end{equation}

Next, we assume that $V$ is quasi-periodic and compute the mean
value of $g(z,\cdot)^{-1}$ using \eqref{3.8}. Before embarking on this
task we briefly review a few properties of quasi-periodic functions. 

We denote by $CP(\bbR)$ and $QP(\bbR)$, the sets of continuous periodic
and quasi-periodic functions on $\bbR$, respectively. In particular, $f$ is
called quasi-periodic with fundamental periods $(\Omega_1,\dots,\Omega_N) 
\in (0,\infty)^N$ if the frequencies 
$2\pi/\Omega_1,\dots,2\pi/\Omega_N$ are linearly independent over $\bbQ$
and if there exists a continuous function
$F\in C(\bbR^N)$, periodic of period $1$ in each of its arguments 
\begin{equation}
F(x_1,\dots,x_j+1,\dots,x_N)=F(x_1,\dots,x_N), \quad
x_j\in\bbR, \;  j=1,\dots,N, \lb{3.14a}
\end{equation}
such that
\begin{equation}
f(x)=F(\Omega_1^{-1}x,\dots,\Omega_N^{-1}x), \quad x\in\bbR. \lb{3.14b}
\end{equation} 
The frequency module $\Mod\,(f)$ of $f$ is then of the type 
\begin{equation}
\Mod\,(f)=\{2\pi m_1/\Omega_1+\cdots +2\pi m_N/\Omega_N\,|\, m_j\in\bbZ,
\, j=1,\dots,N\}. \lb{3.14c}
\end{equation} 

We note that $f\in CP(\bbR)$ if and only if there are
$r_j\in\bbQ\backslash\{0\}$ such that $\Omega_j=r_j\hatt \Omega$ for some
$\hatt \Omega>0$, or equivalently, if and only if $\Omega_j=m_j\wti
\Omega$, $m_j\in\bbZ\backslash\{0\}$ for some $\wti \Omega>0$. $f$ has the
fundamental period $\Omega>0$ if every period of $f$ is an integer
multiple of $\Omega$.

For any quasi-periodic (in fact, Bohr (uniformly) almost periodic)
function $f$, the mean value $\langle f\rangle$ of $f$, defined by
\begin{equation}
\langle f\rangle =\lim_{R\to\infty}\f{1}{2R} \int_{x_0-R}^{x_0+R} dx\,
f(x),
\lb{3.15}
\end{equation}
exists and is independent of $x_0\in\bbR$. Moreover, we recall the
following facts (also valid for Bohr (uniformly) almost periodic
functions on $\bbR$), see, for instance, \cite[Ch.\ I]{Be54},
\cite[Sects.\ 39--92]{Bo47}, \cite[Ch.\ I]{Co89}, \cite[Chs.\ 1,3,6]{Fi74},
\cite{JM82}, \cite[Chs.\ 1,2,6]{LZ82}, and \cite{Sc65}.

\begin{theorem} \lb{t3.3}
Assume $f,g \in QP(\bbR)$ and $x_0,x\in\bbR$. Then the following assertions
hold: \\
$(i)$ $f$ is uniformly continuous on $\bbR$ and $f\in
L^\infty(\bbR;dx)$.\\
$(ii)$ $\ol f$, $d\,f$, $d\in\bbC$, $f(\cdot+c)$, $f(c\cdot)$, $c\in\bbR$,
$|f|^\alpha$, $\alpha\geq 0$ are all in $QP(\bbR)$. \\
$(iii)$ $f+g, fg\in QP(\bbR)$. \\
$(iv)$ $f/g\in QP(\bbR)$ if and only if $\inf_{s\in\bbR}[|g(s)|]>0$. \\
$(v)$ Let $G$ be uniformly continuous on $\calM\subseteq \bbR$ and
$f(s)\in\calM$ for all $s\in\bbR$. Then $G(f)\in QP(\bbR)$. \\
$(vi)$ $f'\in QP(\bbR)$ if and only if $f'$ is uniformly continuous
on $\bbR$. \\
$(vii)$ Let $\langle f\rangle=0$, then $\int_{x_0}^x dx'\,
f(x')\underset{|x|\to\infty}{=}\oh(|x|)$.\\
$(viii)$ Let $F(x)=\int_{x_0}^x dx'\, f(x')$. Then 
$F\in QP(\bbR)$ if and only if $F\in L^\infty(\bbR;dx)$. \\
$(ix)$ If $0\leq f\in QP(\bbR)$, $f\not\equiv 0$, then $\langle
f\rangle>0$. \\
$(x)$ If $f=|f|\exp(i\varphi)$, then $|f|\in QP(\bbR)$ and $\varphi$ is of
the type $\varphi(x)=cx+\psi(x)$, where $c\in \bbR$ and $\psi\in QP(\bbR)$
$($and real-valued\,$)$. \\
$(xi)$ If $F(x)=\exp\Big(\int_{x_0}^{x} dx'\, f(x')\Big)$, then
$F\in QP(\bbR)$ if and only if $f(x)=i\beta + \psi(x)$, where
$\beta\in\bbR$, $\psi\in QP(\bbR)$, and $\Psi\in L^\infty(\bbR;dx)$, where
$\Psi(x)=\int_{x_0}^{x} dx' \, \psi(x')$. 
\end{theorem}

For the rest of this section and the next it will be convenient to
introduce the following hypothesis:

\begin{hypothesis} \lb{h3.3}
Assume the affine part of $\calK_n$ to be nonsingular. Moreover, suppose
that $V\in C^\infty(\bbR)\cap QP(\bbR)$ satisfies the $n$th stationary KdV
equation \eqref{1.2.7} on $\bbR$. 
\end{hypothesis}

Next, we note the following result.

\begin{lemma} \lb{l3.3}
Assume Hypothesis \ref{h3.3}. Then $V^{(k)}$, $k\in\bbN$, and $f_\ell$,
$\ell\in\bbN$, and hence all $x$ and $z$-derivatives of
$F_n(z,\cdot)$, $z\in\bbC$, and $g(z,\cdot)$, $z\in\Pi$, are
quasi-periodic. Moreover, taking limits to points on
$\calC$, the last result extends to either side of the cuts in the set 
$\calC\backslash\{E_m\}_{m=0}^{2n}$ $($cf.\ \eqref{b2}$)$ by 
continuity with respect to $z$.
\end{lemma}
\begin{proof}
Since by hypothesis $V\in C^\infty(\bbR)\cap L^\infty(\bbR;dx)$,
$\sKdV_n(V)=0$ implies $V^{(k)}\in L^\infty(\bbR;dx)$, $k\in\bbN$ and 
$f_\ell\in C^\infty(\bbR)\cap L^\infty(\bbR;dx)$, $\ell\in\bbN_0$,
applying  Remark \ref{r2.7}. In particular $V^{(k)}$ is uniformly
continuous on $\bbR$ and hence quasi-periodic for all $k\in\bbN$. Since the
$f_\ell$ are differential polynomials with respect to $V$, also $f_\ell$,
$\ell\in\bbN$ are quasi-periodic. The corresponding assertion for
$F_n(z,\cdot)$ then follows from \eqref{1.2.9} and that for $g(z,\cdot)$
follows from \eqref{3.13}. 
\end{proof} 

For future purposes we introduce the set 
\begin{align}
\Pi_C&=\Pi \Big\backslash \Big\{\{z\in\bbC\,|\, |z|\leq C+1\} \cup 
\{z\in\bbC\,|\, \Re(z)\geq \min_{m=0,\dots,2n}[\Re(E_m)]-1, \no \\ 
& \qquad \qquad \min_{m=0,\dots,2n}[\Im(E_m)]-1\leq \Im(z)\leq
\max_{m=0,\dots,2n}[\Im(E_m)]+1\}\Big\}, \lb{3.16}
\end{align}
where $C>0$ is the constant in \eqref{1.3C}. Moreover, without loss of
generality, we may assume $\Pi_C$ contains no cuts, that is,
\begin{equation}
\Pi_C\cap \calC=\emptyset. 
\end{equation}

\begin{lemma} \lb{l3.4}
Assume Hypothesis \ref{h3.3} and let $z, z_0\in\Pi$. Then
\begin{equation}
\big\langle g(z,\cdot)^{-1}\big\rangle=
-2\int_{z_0}^z dz'\, \langle g(z',\cdot)\rangle + 
\big\langle g(z_0,\cdot)^{-1}\big\rangle, \lb{3.17}
\end{equation}
where the path connecting $z_0$ and $z$ is assumed to lie in the cut
plane $\Pi$. Moreover, by taking limits to points on $\calC$
in \eqref{3.17}, the result \eqref{3.17} extends to either side of the
cuts in the set $\calC$ by continuity with respect to $z$.
\end{lemma}
\begin{proof}
Let $z, z_0\in\Pi_C$. Integrating equation \eqref{3.8} from
$z_0$ to $z$ along a smooth path in $\Pi_C$ yields
\begin{align}
g(z,x)^{-1} - g(z_0,x)^{-1} &= -2\int_{z_0}^z dz'\, g(z',x) +
[g_{xx}(z,x)-g_{xx}(z_0,x)] \no \\
& \quad -\int_{z_0}^z dz'\, \big[g(z',x)^{-1}g_x(z',x)g_z(z',x)\big]_x \no
\\
& = -2\int_{z_0}^z dz'\, g(z',x) +
g_{xx}(z,x)-g_{xx}(z_0,x) \no \\
& \quad -\bigg[\int_{z_0}^z dz'\, g(z',x)^{-1}g_x(z',x)g_z(z',x)\bigg]_x.
\lb{3.18}
\end{align}
By Lemma \ref{l3.3} $g(z,\cdot)$ and all its $x$-derivatives are
quasi-periodic, 
\begin{equation}
\langle g_{xx}(z,\cdot)\rangle=0, \quad z\in\Pi. \lb{3.19}
\end{equation}
Since we actually assumed $z\in\Pi_C$, also $g(z,\cdot)^{-1}$ 
is quasi-periodic. Consequently, also
\begin{equation}
\int_{z_0}^z dz'\, g(z',\cdot)^{-1}g_x(z',\cdot)g_z(z',\cdot), \quad 
z\in\Pi_C, 
\end{equation}
is a family of uniformly almost periodic functions for $z$
varying in compact subsets of $\Pi_C$ as discussed in 
\cite[Sect.\ 2.7]{Fi74} and one obtains
\begin{equation}
\bigg\langle \bigg[\int_{z_0}^z dz'\,
g(z',\cdot)^{-1}g_x(z',\cdot)g_z(z',\cdot)\bigg]_x\bigg\rangle =0.
\lb{3.21}
\end{equation}
Hence, taking mean values in \eqref{3.18} (taking into account \eqref{3.19}
and \eqref{3.21}), proves \eqref{3.17} for $z\in\Pi_C$. Since $f_\ell$,
$\ell\in\bbN_0$, are quasi-periodic by Lemma \ref{l3.3} (we recall that
$f_0=1$), \eqref{1.2.9} and \eqref{3.12} yield  
\begin{equation}
\int_{z_0}^z dz'\, \langle g(z',\cdot)\rangle = \f{i}{2}
\sum_{\ell=0}^n\langle f_{n-\ell} \rangle \int_{z_0}^z dz'\,
\f{{(z')}^\ell}{R_{2n+1}(z')^{1/2}}.  \lb{3.22}
\end{equation}
Thus, $\int_{z_0}^z dz'\, \langle g(z',\cdot)\rangle$ has an analytic
continuation with respect to $z$ to all of $\Pi$ and consequently,
\eqref{3.17} for $z\in\Pi_C$ extends by analytic continuation to
$z\in\Pi$. By continuity this extends to either side of the cuts in
$\calC$. Interchanging the role of $z$ and $z_0$, analytic
continuation with respect to $z_0$ then yields \eqref{3.17} for 
$z,z_0\in\Pi$.
\end{proof} 

\begin{remark} \lb{r3.7}
For $z\in\Pi_C$, $g(z,\cdot)^{-1}$ is quasi-periodic and hence
$\big\langle g(z,\cdot)^{-1}\big\rangle$ is well-defined. If one
analytically continues $g(z,x)$ with respect to $z$, $g(z,x)$ will acquire
zeros for some $x\in\bbR$ and hence $g(z,\cdot)^{-1}\notin QP(\bbR)$.
Nevertheless, as shown by the right-hand side of \eqref{3.17},
$\big\langle g(z,\cdot)^{-1}\big\rangle$ admits an analytic continuation
in $z$ from $\Pi_C$ to all of $\Pi$, and from now on,
$\big\langle g(z,\cdot)^{-1}\big\rangle$, $z\in\Pi$, always denotes that
analytic continuation (cf.\ also \eqref{3.23}).
\end{remark}

Next, we will invoke the Baker--Akhiezer function $\psi(P,x,x_0)$
and analyze the expression $\big\langle g(z,\cdot)^{-1}\big\rangle$ in
more detail. 

\begin{theorem} \lb{t3.7}
Assume Hypothesis \ref{h3.3}, let $P=(z,y)\in \Pi_\pm$, and
$x,x_0\in\bbR$. Moreover, select a homology basis 
$\{\ti a_j, \ti b_j\}_{j=1}^n$ on $\calK_n$ such that 
$\wti{\ul B}=i\wti{\ul U}^{(2)}_0$, with $\wti{\ul U}^{(2)}_0$ the
vector of $\ti b$-periods of the normalized differential of the second
kind, $\wti \ome_{P_\infty,0}^{(2)}$, satisfies the constraint
\begin{equation}
\wti{\ul B}=i \wti{\ul U}^{(2)}_0 \in \bbR^n 
\end{equation}
$($cf.\ Appendix \ref{sB}$)$. Then,
\begin{equation}
\Re\big(\big\langle g(P,\cdot)^{-1}\big\rangle\big)
= -2\Im\big(y\big\langle F_n(z,\cdot)^{-1}\big\rangle\big)
=2\Im\bigg(\int_{Q_0}^P \wti \ome_{P_\infty,0}^{(2)} 
-\ti e^{(2)}_0(Q_0)\bigg).  \lb{3.23}
\end{equation}
\end{theorem}
\begin{proof}
Using \eqref{1.3.20}, one obtains for $z\in\Pi_C$,
\begin{align}
\psi (P,x,x_0) &=\bigg(\f{F_n (z,x)}{F_n (z,x_0)} \bigg)^{1/2} \exp
\bigg(iy \int_{x_0}^x  dx' \, F_n (z,x')^{-1} \bigg) \no \\
&=\bigg(\f{F_n (z,x)}{F_n (z,x_0)} \bigg)^{1/2} \exp
\bigg(iy \int_{x_0}^x  dx' \, \big[F_n (z,x')^{-1} - 
\big\langle F_n (z,\cdot)^{-1}\big\rangle \big]\bigg) \no \\
& \quad \times \exp\big(i(x-x_0)y\big\langle F_n
(z,\cdot)^{-1}\big\rangle\big), \lb{3.26} \\
& \quad P=(z,y)\in\Pi_\pm, \; z\in\Pi_C, \; x,x_0\in\bbR. \no
\end{align}
Since $\big[F_n (z,x')^{-1} - \big\langle F_n (z,\cdot)^{-1}\big\rangle
\big]$ has mean zero,
\begin{equation}
\bigg|\int_{x_0}^x  dx' \, \big[F_n (z,x')^{-1} - 
\langle F_n (z,\cdot)^{-1}\rangle \big]\bigg|\underset{|x|\to\infty}{=} 
\oh(|x|), \quad z\in\Pi_C  \lb{3.27} 
\end{equation}
by Theorem \ref{t3.3}\,(vii). In addition, the factor $F_n (z,x)/F_n
(z,x_0)$ in
\eqref{3.26} is quasi-periodic and hence bounded on $\bbR$. 

On the other hand, \eqref{1.3.59} yields
\begin{align}
\psi(P,x,x_0) &= \frac{\theta (\uz(P_\infty,\humu (x_0)))\theta
(\uz(P,\humu (x)))}{\theta (\uz(P_\infty,\humu (x)))\theta (\uz(P,\humu
(x_0)))} \no \\
& \quad \times \exp \bigg[-i (x-x_0)\bigg(\int_{Q_0}^P
\wti \ome_{P_\infty,0}^{(2)}-\ti e^{(0)}_0(Q_0)\bigg)\bigg] \no \\
&= \Theta(P,x,x_0)\exp \bigg[-i (x-x_0)\bigg(\int_{Q_0}^P
\wti\ome_{P_\infty,0}^{(2)}-\ti e^{(2)}_0(Q_0)\bigg)\bigg], \lb{3.28}
\\ 
& \hspace*{3.9cm} 
P\in\calK_n\backslash\big\{\{\Pinf\}\cup\{\hmu_j(x_0)\}_{j=1}^n\big\}.
\no 
\end{align}
Taking into account \eqref{1.3.60d}, \eqref{1.3.60}, \eqref{2.70}, 
\eqref{b9}, and the fact that by \eqref{1.3C} no $\hmu_j(x)$ can reach
$\Pinf$ as $x$ varies in $\bbR$, one concludes that 
\begin{equation}
\Theta(P,\cdot,x_0) \in L^\infty(\bbR;dx), \quad 
P\in\calK_n\backslash \{\hmu_j(x_0)\}_{j=1}^n.
\lb{3.29}
\end{equation}
A comparison of \eqref{3.26} and \eqref{3.28} then shows that the
$\oh(|x|)$-term in \eqref{3.27} must actually be bounded on $\bbR$ and
hence the left-hand side of \eqref{3.27} is quasi-periodic. In
addition, the term
\begin{equation}
\exp\bigg(iR_{2n+1}(z)^{1/2} \int_{x_0}^x  dx' \, \big[F_n (z,x')^{-1} - 
\big\langle F_n (z,\cdot)^{-1}\big\rangle \big]\bigg), \quad 
z\in\Pi_C, \lb{3.31}
\end{equation}
is then quasi-periodic by Theorem \ref{t3.3}\,$(xi)$. A further comparison
of \eqref{3.26} and \eqref{3.28} then yields \eqref{3.23} for $z\in\Pi_C$.
Analytic continuation with respect to $z$ then yields \eqref{3.23} for
$z\in\Pi$. By continuity  with respect to $z$, taking 
boundary values to either side  of the cuts in the set $\calC$, this then
extends to $z\in\calC$ (cf.\ \eqref{b2}, \eqref{b3}) and hence proves
\eqref{3.23} for $P=(z,y)\in \calK_n\backslash \{ P_\infty\}$. 
\end{proof}

\section{Spectra of Schr\"odinger operators with
quasi-periodic algebro-geometric KdV  potentials} \lb{s4}

In this section we establish the connection between the 
algebro-geometric formalism of Section \ref{s2} and the spectral
theoretic description of Schr\"odinger operators $H$ in
$L^2(\bbR;dx)$ with quasi-periodic algebro-geometric KdV potentials.
In particular, we introduce the conditional stability set of $H$
and prove our principal result, the characterization of the spectrum of
$H$. Finally, we provide a qualitative description of the spectrum of
$H$ in terms of analytic spectral arcs.

Suppose that $V\in C^\infty(\bbR)\cap QP(\bbR)$ satisfies the $n$th
stationary KdV equation \eqref{1.2.7} on $\bbR$. The corresponding
Schr\"odinger operator $H$ in $L^2(\bbR;dx)$ is then introduced by
\begin{equation}
H=-\f{d^2}{dx^2}+V, \quad \dom(H)=H^{2,2}(\bbR). \lb{3.35}
\end{equation}
Thus, $H$ is a densely defined closed operator in $L^2(\bbR;dx)$ (it is
self-adjoint if and only if $V$ is real-valued).

Before we turn to the spectrum of $H$ in the general non-self-adjoint
case, we briefly mention the following result on the spectrum of $H$
in the self-adjoint case with a quasi-periodic (or almost
periodic) real-valued potential $q$. We denote by $\sigma(A)$,
$\sigma_{\rm{e}}(A)$, and $\sigma_{\rm d}(A)$ the spectrum, essential
spectrum, and discrete spectrum of a self-adjoint operator $A$ in a
complex Hilbert space, respectively.

\begin{theorem} [See, e.g., \cite{Si82}] 
Let $V\in QP(\bbR)$ and $q$ be real-valued. Define the self-adjoint
Schr\"odinger operator $H$ in $L^2(\bbR;dx)$ as in \eqref{3.35}. Then, 
\begin{equation}
\sigma(H)=\sigma_{\rm{e}}(H)\subseteq \big[\min_{x\in\bbR}
(V(x)),\infty\big), \quad \sigma_{\rm d}(H)=\emptyset. \lb{3.15b}
\end{equation}
Moreover, $\sigma(H)$ contains no isolated points, that is, $\sigma(H)$
is a perfect set. 
\end{theorem}

In the special periodic case where $V\in CP(\bbR)$ is real-valued, the
spectrum  of $H$ is purely absolutely continuous and either a finite
union of some compact intervals and a half-line or an infinite union
of compact intervals (see, e.g., \cite[Sect.\ 5.3]{Ea73}, \cite[Sect.\
XIII.16]{RS78}). If $V\in CP(\bbR)$ and $V$ is complex-valued, then the
spectrum of $H$ is purely continuous and it consists of either a
finite union of simple analytic arcs and one simple semi-infinite
analytic arc tending to infinity or an infinite union of simple
analytic arcs (cf.\ \cite{Ro63}, \cite{Se60}, and 
\cite{Tk64})\footnote{in either case the resolvent set is 
connected.}. 

\begin{remark} \lb{r4.2}
Here $\sigma\subset\bbC$ is called an {\it arc} if there exists a
parameterization $\gamma\in C([0,1])$ such that
$\sigma=\{\gamma(t)\,|\, t\in [0,1]\}$. The arc $\sigma$ is called
{\it simple} if there exists a parameterization $\gamma$
such that $\gamma\colon [0,1]\to\bbC$ is injective. The arc $\sigma$
is called {\it analytic} if there is a parameterization $\gamma$ that
is analytic at each $t\in [0,1]$. Finally, $\sigma_\infty$ is called a
{\it semi-infinite} arc if there exists a parameterization $\gamma\in
C([0,\infty))$ such that 
$\sigma_\infty=\{\gamma(t)\,|\, t\in [0,\infty)\}$ and $\sigma_\infty$
is an unbounded subset of $\bbC$. Analytic semi-infinite arcs are
defined analogously and by a simple semi-infinite arc we mean one that
is without self-intersection (i.e., corresponds to a injective
parameterization) with the additional restriction that the
unbounded part of $\sigma_\infty$ consists of precisely one branch
tending to infinity.
\end{remark}

Now we turn to the analyis of the generally non-self-adjoint operator
$H$ in \eqref{3.35}. Assuming Hypothesis \ref{h3.3} we now introduce
the set $\Sigma\subset\bbC$ by
\begin{equation}
\Sigma=\big\{\lambda\in\bbC\,\big|\, \Re\big(\big\langle
g(\lambda,\cdot)^{-1}\big\rangle\big)=0\big\}.
\lb{3.34}
\end{equation}
Below we will show that $\Sigma$ plays the role of the
conditional stability set of $H$, familiar from the spectral theory of
one-dimensional periodic Schr\"odinger operators (cf.\ 
\cite[Sect.\ 5.3]{Ea73}, \cite{Ro63}, \cite{We98}, \cite{We98a}).

\begin{lemma} \lb{l.3.9}
Assume Hypothesis \ref{h3.3}. Then $\Sigma$ coincides with the
conditional stability set of $H$, that is,
\begin{align}
\Sigma&=\{\lambda\in\bbC\,|\, \text{there exists at least one bounded
distributional solution} \no \\
& \hspace*{1.8cm} \text{$0\neq\psi\in L^\infty(\bbR;dx)$ of
$H\psi=\lambda\psi$.}\}  \lb{3.36}
\end{align} 
\end{lemma}
\begin{proof}
By \eqref{3.28} and \eqref{3.29}, 
\begin{align}
\psi(P,x)&=
\frac{\theta(\uz(P,\humu (x)))}{\theta (\uz(P_\infty,\humu (x)))} 
 \exp \bigg[-i x\bigg(\int_{Q_0}^P \wti\ome_{P_\infty,0}^{(2)} 
-\ti e^{(0)}_0(Q_0)\bigg)\bigg], \lb{3.36a} \\ 
& \hspace*{5.9cm} P=(z,y)\in\Pi_\pm, \no
\end{align}
is a distributional solution of $H\psi=z\psi$ which is bounded on $\bbR$
if and only if the exponential function in \eqref{3.36a} is bounded on
$\bbR$. By \eqref{3.23}, the latter holds if and only if 
\begin{equation}
\Re\big(\big\langle
g(z,\cdot)^{-1}\big\rangle\big)=0. \lb{3.36b}
\end{equation} 
\end{proof}

\begin{remark} \lb{r3.12} 
At first sight our {\it a priori} choice of cuts $\calC$ for
$R_{2n+1}(\cdot)^{1/2}$, as described in Appendix \ref{sA}, might seem
unnatural as they completely ignore the actual spectrum of $H$. However,
the spectrum of $H$ is not known from the outset, and in the case of
complex-valued periodic potentials, spectral arcs of $H$ may actually
cross each other (cf.\ \cite{GW95}, \cite{PT91}, and 
Theorem \ref{t3.14}\,(iv)) which renders them unsuitable for cuts of
$R_{2n+1}(\cdot)^{1/2}$. 
\end{remark}

Before we state our first principal result on the spectrum of $H$, we 
find it convenient to recall a number of basic definitions and well-known
facts in connection with the spectral theory of  non-self-adjoint
operators (we refer to \cite[Chs.\ I, III, IX]{EE89},  
\cite[Sects.\ 1, 21--23]{Gl65}, \cite[Sects.\ IV.5.6, V.3.2]{Ka80}, and
\cite[p.\ 178--179]{RS78} for more details). Let $S$ be a densely defined
closed operator in a complex separable Hilbert space $\cH$. Denote by
$\cB(\cH)$ the Banach space of all bounded linear operators on $\cH$ and
by $\ker(T)$ and $\ran(T)$ the kernel (null space) and range of a linear
operator $T$ in $\cH$. The resolvent set,
$\rho(S)$, spectrum, $\sigma(S)$, point spectrum (the set of
eigenvalues), 
$\sigma_{\rm p}(S)$, continuous spectrum, $\sigma_{\rm c}(S)$, residual
spectrum, $\sigma_{\rm r}(S)$, field of regularity, $\pi(S)$, approximate 
point spectrum, $\sigma_{\rm ap}(S)$, two kinds of essential spectra,
$\sigma_{\rm e}(S)$, and $\wti\sigma_{\rm e}(S)$, the numerical range
of $S$, $\Theta(S)$, and the sets $\Delta (S)$ and $\wti\Delta (S)$ are
defined as follows: 
\begin{align}
\rho(S)&=\{z\in\bbC\,|\, (S-z I)^{-1}\in \cB(\cH)\},
\lb{5.15} \\
\sigma(S)&=\bbC\backslash\rho (S), \lb{5.8} \\
\sigma_{\rm p}(S)&=\{\lambda\in\bbC\,|\, \ker(S-\lambda I)\neq\{0\} \}, 
\lb{5.9} \\
\sigma_{\rm c}(S)&=\{\lambda\in\bbC \,|\, \text{$\ker(S-\lambda
I)=\{0\}$ and $\ran(S-\lambda I)$ is dense in $\cH$} \no \\ &
\hspace*{1.8cm} \text{but not equal to
$\cH$}\}, \lb{5.10} \\
\sigma_{\rm r}(S)&=\{\lambda\in\bbC\,|\, \text{$\ker(S-\lambda I)=\{0\}$  
and $\ran(S-\lambda I)$ is not dense in $\cH$}\}, \lb{5.11} \\
\pi(S)&=\{z \in\bbC \,|\, \text{there exists $k_z >0$ s.t. 
$\| (S- zI)u\|_\cH \ge k_z \| u\|_\cH$} \no \\
& \hspace*{1.75cm} \text{for all $u\in\dom(S)$}\}, \lb{5.14} \\
\sigma_{\rm ap}(S)&=\bbC\backslash\pi(S), \lb{5.21} \\
\Delta(S)&=\{z\in\bbC \,|\, \text{$\dim(\ker(S-zI))<\infty$ and
$\ran(S-zI)$ is closed}\}, \lb{5.16} \\
\sigma_{\rm e}(S)&= \bbC\backslash\Delta (S), \lb{5.22b} \\
\wti\Delta(S)&=\{z\in\bbC \,|\, \text{$\dim(\ker(S-zI))<\infty$ or 
$\dim(\ker(S^*-\ol z I))<\infty$}\}, \lb{5.16a} \\
\wti\sigma_{\rm e} (S)&=\bbC\backslash \wti\Delta(S), \lb{5.16b} \\
\Theta(S)&=\{(f,Sf)\in\bbC \,|\, f\in\dom(S), \, \|f\|_{\cH}=1\},
\lb{5.16c}
\end{align}
respectively. One then has 
\begin{align}
\sigma (S)&=\sigma_{\rm p}(S)\cup\sigma_{\rm{c}}(S)\cup
\sigma_{\rm r}(S) \quad \text{(disjoint union)} \lb{5.17} \\
&=\sigma_{\rm p}(S)\cup\sigma_{\rm{e}}(S)\cup\sigma_{\rm r}(S), 
\lb{5.18} \\
\sigma_{\rm c}(S)&\subseteq\sigma_{\rm e}(S)\backslash
(\sigma_{\rm p}(S)\cup\sigma_{\rm r}(S)), \lb{5.18a} \\
\sigma_{\rm r}(S)&=\sigma_{\rm p}(S^*)^* \backslash\sigma_{\rm p}(S), 
\lb{5.19} \\
\sigma_{\rm ap}(S)&=\{\lambda \in \bbC \, |\, \text{there exists a
sequence $\{ f_n\}_{n\in\bbN}\subset\dom(S)$} \no \\ 
&\hspace*{1.75cm} \text{with $\| f_n \|_\cH=1$, $n\in\bbN$, and
$\lim_{n\to\infty} \|(S-\lambda I)f_n\|_\cH=0$}\}, \lb{5.12} \\
\wti\sigma_{\rm e}(S)&\subseteq \sigma_{\rm e}(S)\subseteq\sigma_{\rm
ap}(S)\subseteq\sigma(S) \, \text{ (all four sets are closed)}, \lb{5.23}
\\
\rho(S)&\subseteq \pi(S) \subseteq \Delta (S) \subseteq \wti\Delta (S) 
\;\; \text{ (all four sets are open),} \lb{5.24} \\
\wti\sigma_{\rm e}(S) & \subseteq \ol{\Theta(S)}, \quad \Theta(S) \,
\text{ is convex,} \lb{5.25} \\
\wti\sigma_{\rm e}(S) &=\sigma_{\rm e}(S) \, \text{ if $S=S^*$.} \lb{5.26}
\end{align}
Here $\sigma^*$ in the context of \eqref{5.19} denotes the complex 
conjugate of the set $\sigma\subseteq\bbC$, that is,
\begin{equation}
\sigma^*=\{\ol{\lambda}\in\bbC \,|\, \lambda\in\sigma\}. \lb{5.27}
\end{equation}
We note that there are several other versions of the concept of the
essential spectrum in the non-self-adjoint context (cf.\
\cite[Ch.\ IX]{EE89}) but we will only use the two in \eqref{5.22b} 
and in \eqref{5.16b} in this paper.

Finally, we recall the following result due to Talenti \cite{Ta69}
and Tomaselli \cite{To69} (see also Chisholm and Everitt
\cite{CE70}, Chisholm, Everitt, and Littlejohn \cite{CEL99}, and
Muckenhoupt \cite{Mu72}).

\begin{lemma}  \lb{l3.10} 
Let $f\in L^2(\bbR;dx)$, $U\in L^2((-\infty,R];dx)$, and $V\in
L^2([R,\infty);dx)$ for all $R\in\bbR$. Then the following assertions
$(i)$--$(iii)$ are equivalent: \\
$(i)$ There exists a finite constant $C>0$ such that
\begin{equation}
\int_\bbR dx \, \bigg|U(x)\int_x^\infty dx' \,V(x')f(x')
\bigg|^2 \leq C \int_\bbR dx\, |f(x)|^2. \lb{3.37}
\end{equation}
$(ii)$ There exists a finite constant $D>0$ such that
\begin{equation}
\int_\bbR dx \, \bigg|V(x)\int_{-\infty}^x dx' \,U(x')f(x')
\bigg|^2 \leq D \int_\bbR dx\, |f(x)|^2. \lb{3.38}
\end{equation}
$(iii)$ 
\begin{equation}
\sup_{r\in\bbR}\Bigg[\bigg(\int_{-\infty}^r dx\, |U(x)|^2\bigg)
\bigg(\int_r^\infty dx \, |V(x)|^2\bigg)\Bigg] <\infty. \lb{3.39}
\end{equation}
\end{lemma}

We start with the following elementary result.

\begin{lemma}  \lb{l3.10a} 
Let $H$ be defined as in \eqref{3.35}. Then,
\begin{equation}
\sigma_{\rm e}(H)=\wti\sigma_{\rm e}(H)\subseteq \ol{\Theta(H)}.
\lb{3.39a}
\end{equation}
\end{lemma}
\begin{proof}
Since $H$ and $H^*$ are second-order ordinary differential
operators on $\bbR$, 
\begin{equation}
\dim(\ker(H-z I))\leq 2, \quad \dim(\ker(H^*-\ol z I))\leq 2.
\lb{3.39b}
\end{equation}
Equations \eqref{5.16}--\eqref{5.16b} and \eqref{5.25} then prove
\eqref{3.39a}.
\end{proof}

\begin{theorem}  \lb{t3.11} 
Assume Hypothesis \ref{h3.3}. Then the point spectrum
and residual spectrum of $H$ are empty and hence the spectrum of $H$ is
purely continuous,
\begin{align}
&\sigma_{\rm p}(H)=\sigma_{\rm r}(H)=\emptyset,  \lb{3.40} \\
&\sigma(H)=\sigma_{\rm c}(H)=\sigma_{\rm e}(H)=\sigma_{\rm ap}(H).
\lb{3.41}
\end{align}
\end{theorem}
\begin{proof}
First we prove the absence of the point spectrum of $H$. Suppose
$z\in\Pi\backslash\{\Sigma \cup\{\mu_j(x_0)\}_{j=1}^n\}$.
Then $\psi(P,\cdot,x_0)$ and $\psi(P^*,\cdot,x_0)$ are linearly
independent distributional solutions of $H\psi=z\psi$ which are unbounded
at $+\infty$ or $-\infty$. This argument extends to all
$z\in\Pi\backslash\Sigma$ by multiplying $\psi(P,\cdot,x_0)$ and
$\psi(P^*,\cdot,x_0)$ with an appropriate function of $z$ and $x_0$
(independent of $x$). It also extends to either side of the cut
$\calC\backslash\Sigma$ by continuity with respect to $z$. On the
other hand, since $V^{(k)}\in L^\infty(\bbR;dx)$ for all $k\in\bbN_0$, any
distributional solution $\psi(z,\cdot)\in L^2(\bbR;dx)$ of $H\psi=z\psi$,
$z\in\bbC$, is necessarily bounded. In fact, 
\begin{equation}
\psi^{(k)}(z,\cdot)\in L^\infty(\bbR;dx)\cap L^2(\bbR;dx), 
\quad k\in\bbN_0, \lb{3.68A}
\end{equation} 
applying $\psi''(z,x)=(V(x)-z)\psi(z,x)$ and \eqref{2.58} with $p=2$ and
$p=\infty$ repeatedly. (Indeed, $\psi(z,\cdot)\in L^2(\bbR;dx)$ implies 
$\psi''(z,\cdot)\in L^2(\bbR;dx)$ which in turn implies $\psi'(z,\cdot)\in
L^2(\bbR;dx)$. Integrating $(\psi^2)'=2\psi\psi'$ then yields 
$\psi(z,\cdot)\in L^\infty (\bbR;dx)$. The latter yields 
$\psi''(z,\cdot)\in L^\infty (\bbR;dx)$, etc.) Thus, 
\begin{equation}
\{\bbC\backslash\Sigma\}\cap \sigma_{\rm p}(H)=\emptyset. \lb{3.67a}
\end{equation}
Hence, it remains to rule out eigenvalues located in $\Sigma$. 
We consider a fixed $\lambda\in\Sigma$ and note that by \eqref{1.3.19},
there exists at least one distributional solution $\psi_1(\lambda,\cdot)
\in L^\infty(\bbR;dx)$ of
$H\psi=\lambda\psi$. Actually, a comparison of \eqref{1.3.20} and
\eqref{3.34} shows that we may choose $\psi_1(\lambda,\cdot)$ such that
$|\psi_1(\lambda,\cdot)|\in QP(\bbR)$ and hence
$\psi_1(\lambda,\cdot)\notin L^2(\bbR;dx)$. As in \eqref{3.68A} one then
infers from repeated use of $\psi''(\lambda)=(V-\lambda)\psi(\lambda)$ and
\eqref{2.58} with $p=\infty$ that 
\begin{equation}
\psi_1^{(k)}(\lambda,\cdot)\in L^\infty(\bbR;dx), \quad k\in\bbN_0.  
\lb{3.68a}
\end{equation}
Next, suppose there exists a second distributional solution
$\psi_2(\lambda,\cdot)$ of $H\psi=\lambda\psi$ which is linearly
independent of $\psi_1(\lambda,\cdot)$ and which satisfies 
$\psi_2(\lambda,\cdot)\in L^2(\bbR;dx)$. Applying \eqref{3.68A} then
yields
\begin{equation}
\psi_2^{(k)}(\lambda,\cdot)\in L^2(\bbR;dx), \quad k\in\bbN_0. \lb{3.68b}
\end{equation}
Combining \eqref{3.68a} and \eqref{3.68b}, one concludes that the
Wronskian of $\psi_1(\lambda,\cdot)$ and $\psi_2(\lambda,\cdot)$ lies in
$L^2(\bbR;dx)$,
\begin{equation}
W(\psi_1(\lambda,\cdot),\psi_2(\lambda,\cdot))\in L^2(\bbR;dx). \lb{3.68c}
\end{equation}
However, by hypothesis, 
$W(\psi_1(\lambda,\cdot),\psi_2(\lambda,\cdot))=c(\lambda)\neq 0$
is a nonzero constant. This contradiction proves that 
\begin{equation}
\Sigma\cap \sigma_{\rm p}(H)=\emptyset \lb{3.68d}
\end{equation}
and hence $\sigma_{\rm p}(H)=\emptyset$.

Next, we note that the same argument yields that $H^*$ also has no point
spectrum,
\begin{equation}
\sigma_{\rm p}(H^*)=\emptyset. \lb{3.43}
\end{equation}
Indeed, if $V\in C^\infty(\bbR)\cap QP(\bbR)$ satisfies the $n$th
stationary KdV equation \eqref{1.2.7} on $\bbR$, then $\ol V$ also 
satisfies one of the $n$th stationary KdV equations \eqref{1.2.7}
associated with a hyperelliptic curve of genus $n$ with
$\{E_m\}_{m=0}^{2n}$ replaced by $\{\ol E_m\}_{m=0}^{2n}$, etc.
Since by general principles (cf.\ \eqref{5.27}), 
\begin{equation}
\sigma_{\rm r}(B)\subseteq \sigma_{\rm p}(B^*)^* \lb{3.44}
\end{equation}
for any densely defined closed linear operator $B$ in some complex
separable Hilbert space (see, e.g., \cite[p.\ 71]{Go85}), one obtains 
$\sigma_{\rm r}(H)=\emptyset$ and hence \eqref{3.40}. This proves that
the spectrum of $H$ is purely continuous, $\sigma(H)=\sigma_{\rm c}(H)$. 
The remaining equalities in \eqref{3.41} then follow from \eqref{5.18a}
and \eqref{5.23}. 
\end{proof}

The following result is a fundamental one:

\begin{theorem}  \lb{t3.12} 
Assume Hypothesis \ref{h3.3}. Then the spectrum of $H$ coincides with
$\Sigma$ and hence equals the conditional stability set of $H$, 
\begin{align}
\sigma(H) &=\big\{\lambda\in\bbC\,\big|\, \Re\big(\big\langle
g(\lambda,\cdot)^{-1}\big\rangle\big)=0\big\}  \lb{3.45} \\ 
&=\{\lambda\in\bbC\,|\, \text{there exists at least one bounded
distributional solution}  \no \\
& \hspace*{1.8cm} \text{$0\neq\psi\in L^\infty(\bbR;dx)$ of
$H\psi=\lambda\psi$}\}.  \lb{3.45A} 
\end{align}
In particular,
\begin{equation}
\{E_m\}_{m=0}^{2n}\subset\sigma(H), \lb{3.45a}
\end{equation}
and $\sigma(H)$ contains no isolated points.
\end{theorem}
\begin{proof}
First we will prove that 
\begin{equation}
\sigma(H)\subseteq \Sigma \lb{3.42}
\end{equation}
by adapting a method due to Chisholm and Everitt \cite{CE70}. For this
purpose we temporarily choose
$z\in\Pi\backslash\{\Sigma\cup\{\mu_j(x_0)\}_{j=1}^n\}$ and
construct the resolvent of $H$ as follows. Introducing the two
branches $\psi_\pm (P,x,x_0)$ of the Baker--Akhiezer function
$\psi(P,x,x_0)$ by 
\begin{equation}
\psi_\pm (P,x,x_0)=\psi(P,x,x_0), \quad P=(z,y)\in\Pi_\pm, \;
x,x_0\in\bbR, \lb{3.47}
\end{equation} 
we define
\begin{align}
\hat \psi_+(z,x,x_0)&=\begin{cases} \psi_+(z,x,x_0) & \text{if
$\psi_+(z,\cdot,x_0)\in L^2((x_0,\infty);dx)$,} \\
\psi_-(z,x,x_0) & \text{if
$\psi_-(z,\cdot,x_0)\in L^2((x_0,\infty);dx)$,} \end{cases} \lb{3.48} \\
\hat \psi_-(z,x,x_0)&=\begin{cases} \psi_-(z,x,x_0) & \text{if
$\psi_-(z,\cdot,x_0)\in L^2((-\infty,x_0);dx)$,} \\
\psi_+(z,x,x_0) & \text{if
$\psi_+(z,\cdot,x_0)\in L^2((-\infty,x_0);dx)$,} \end{cases} \lb{3.49} \\
& \hspace*{5cm} z\in\Pi\backslash\Sigma, \; x,x_0\in\bbR, \no
\end{align}
and 
\begin{align}
G(z,x,x')&=\f{1}{W(\hat\psi_+(z,x,x_0),\hat\psi_-(z,x,x_0))}\begin{cases}
\hat\psi_-(z,x',x_0)\hat\psi_+(z,x,x_0), & x\geq x', \\
\hat\psi_-(z,x,x_0)\hat\psi_+(z,x',x_0), & x\leq x', \end{cases} \no \\
& \hspace*{6.3cm} z\in\Pi\backslash\Sigma, \; x,x_0\in\bbR. \lb{3.50}
\end{align}
Due to the homogeneous nature of $G$, \eqref{3.50} extends to all
$z\in\Pi$. Moreover, we extend \eqref{3.48}--\eqref{3.50} to either side
of the cut $\calC$ except at possible points in $\Sigma$ (i.e., to 
$\calC\backslash\Sigma$) by continuity with respect to $z$, taking 
limits to $\calC\backslash\Sigma$. Next, we introduce the operator
$R(z)$ in $L^2(\bbR;dx)$ defined by
\begin{align}
(R(z)f)(x)=\int_\bbR dx' \,G(z,x,x')f(x'), \quad f\in C_0^\infty(\bbR), \;
z\in\Pi, \lb{3.51}
\end{align}
and extend it to $z\in\calC\backslash\Sigma$, as
discussed in connection with $G(\cdot,x,x')$. The explicit form of
$\hat\psi_\pm(z,x,x_0)$, inferred from \eqref{3.28} by restricting $P$ to
$\Pi_\pm$, then yields the estimates
\begin{equation} 
|\hat \psi_\pm(z,x,x_0)|\leq C_\pm(z,x_0) e^{\mp \kappa(z)x}, \quad
z\in\Pi\backslash\Sigma, \;  x\in\bbR \lb{3.51a}
\end{equation}
for some constants $C_\pm(z,x_0)>0$, $\kappa(z)>0$,
$z\in\Pi\backslash\Sigma$. An application of Lemma
\ref{l3.10} identifying $U(x)=\exp(-\kappa(z)x)$ and
$V(x)=\exp(\kappa(z)x)$ then proves that $R(z)$,
$z\in\bbC\backslash\Sigma$, extends from $C_0^\infty(\bbR)$ to a bounded
linear operator defined on all of $L^2(\bbR;dx)$. (Alternatively, one can
follow the second part of the proof of Theorem\ 5.3.2 in \cite{Ea73} line
by line.) A straightforward differentiation then proves
\begin{equation}
(H-zI)R(z)f=f, \quad f\in L^2(\bbR;dx), \;
z\in\bbC\backslash\Sigma \lb{3.52}
\end{equation}
and hence also 
\begin{equation}
R(z)(H-zI)g=g, \quad g\in \dom(H), \; z\in\bbC\backslash\Sigma. \lb{3.53}
\end{equation}
Thus, $R(z)=(H-zI)^{-1}$, $z\in\bbC\backslash\Sigma$, and hence
\eqref{3.42} holds.

Next we will prove that 
\begin{equation}
\sigma(H)\supseteq \Sigma. \lb{3.54}
\end{equation}
We will adapt a strategy of proof applied by Eastham in the case of
(real-valued) periodic potentials \cite{Ea67} (reproduced in the proof of
Theorem\ 5.3.2 of \cite{Ea73}) to the (complex-valued) quasi-periodic case
at hand. Suppose $\lambda\in\Sigma$. By the characterization \eqref{3.36}
of $\Sigma$, there exists a bounded distributional solution
$\psi(\lambda,\cdot)$ of $H\psi=\lambda\psi$. A comparison with the
Baker-Akhiezer function \eqref{1.3.20} then shows that we can assume,
without loss of generality, that
\begin{equation}
|\psi(\lambda,\cdot)|\in QP(\bbR). \lb{3.55}
\end{equation} 
Moreover, by the same argument as in the proof of Theorem \ref{t3.11}
(cf.\ \eqref{3.68a}), one obtains
\begin{equation}
\psi^{(k)}(\lambda,\cdot)\in L^\infty(\bbR;dx), \quad k\in\bbN_0.
\lb{3.55a}
\end{equation}
Next, we pick $\Omega>0$ and consider $g\in C^\infty([0,\Omega])$
satisfying
\begin{align}
&g(0)=0, \quad g(\Omega)=1, \no \\
&g'(0)=g''(0)=g'(\Omega)=g''(\Omega)=0, \lb{3.55b} \\
&0\leq g(x)\leq 1, \quad x\in [0,\Omega]. \no
\end{align}
Moreover, we introduce the sequence $\{h_n\}_{n\in\bbN}\in L^2(\bbR;dx)$
by
\begin{equation}
h_n(x)=\begin{cases} 1, & |x|\leq (n-1)\Omega, \\
g(n\Omega-|x|), & (n-1)\Omega\leq |x|\leq n\Omega, \\
0, & |x|\geq n\Omega \end{cases} \lb{3.55c}
\end{equation}
and the sequence $\{f_n(\lambda)\}_{n\in\bbN}\in L^2(\bbR;dx)$ by
\begin{equation}
f_n(\lambda,x)=d_n(\lambda)\psi(\lambda,x)h_n(x), \quad x\in\bbR, \;
d_n(\lambda)>0, \; n\in\bbN. \lb{3.55d}
\end{equation} 
Here $d_n(\lambda)$ is determined by the requirement
\begin{equation}
\|f_n(\lambda)\|_2 =1, \quad n\in\bbN. \lb{3.55e}
\end{equation}
One readily verifies that
\begin{equation}
f_n(\lambda,\cdot) \in\dom(H)=H^{2,2}(\bbR), \quad n\in\bbN. \lb{3.55ea}
\end{equation}
Next, we note that as a consequence of Theorem \ref{t3.3}\,(ix),
\begin{equation}
\int_{-T}^T dx\, |\psi(\lambda,x)|^2\underset{T\to\infty}{=} 2\big\langle
|\psi(\lambda,\cdot)|^2 \big\rangle T +\oh(T) \lb{3.55f}
\end{equation}
with
\begin{equation}
\big\langle |\psi(\lambda,\cdot)|^2 \big\rangle >0. \lb{3.55g}
\end{equation}
Thus, one computes
\begin{align}
1&=\|f_n(\lambda)\|^2_2=d_n(\lambda)^2\int_\bbR dx\, |\psi(\lambda,x)|^2
h_n(x)^2 \no \\
& = d_n(\lambda)^2\int_{|x|\leq n\Omega} dx\, |\psi(\lambda,x)|^2
h_n(x)^2 \geq d_n(\lambda)^2 \int_{|x|\leq (n-1)\Omega} dx\,
|\psi(\lambda,x)|^2 \no \\
& \quad \geq d_n(\lambda)^2 \big[\big\langle |\psi(\lambda,\cdot)|^2
\big\rangle (n-1)\Omega + \oh(n) \big]. \lb{3.55h}
\end{align}
Consequently,
\begin{equation}
d_n(\lambda)\underset{n\to\infty}{=} \Oh\big(n^{-1/2}\big). \lb{3.55i}
\end{equation}
Next, one computes
\begin{equation}
(H-\lambda I)f_n(\lambda,x)=-d_n(\lambda)[2\psi'(\lambda,x)h_n'(x) + 
\psi(\lambda,x)h_n''(x)] \lb{3.55j}
\end{equation}
and hence
\begin{align}
\|(H-\lambda I)f_n\|_2 & \leq d_n(\lambda)[2\|\psi'(\lambda)h_n'\|_2 + 
\|\psi(\lambda)h_n''\|_2], \quad n\in\bbN. \lb{3.55k}
\end{align}
Using \eqref{3.55a} and \eqref{3.55c} one estimates
\begin{align}
\|\psi'(\lambda)h_n'\|_2^2&= \int_{(n-1)\Omega\leq |x|\leq n\Omega} dx \, 
|\psi'(\lambda,x)|^2 |h_n'(x)|^2  
 \leq 2\|\psi'(\lambda)\|_{\infty}^2 \int_0^\Omega dx \, 
|g'(x)|^2 \no \\
&\leq 2 \Omega \|\psi'(\lambda)\|_{\infty}^2 
\|g'\|_{L^\infty([0,\Omega];dx)}^2,  \lb{3.55l}
\end{align}
and similarly,
\begin{align}
\|\psi(\lambda)h_n''\|_2^2&= \int_{(n-1)\Omega\leq |x|\leq n\Omega} dx \, 
|\psi(\lambda,x)|^2 |h_n''(x)|^2 
\leq 2\|\psi(\lambda)\|_{\infty}^2 \int_0^\Omega dx \, 
|g''(x)|^2 \no \\
& \leq 2 \Omega \|\psi(\lambda)\|_{\infty}^2 
\|g''\|_{L^\infty([0,\Omega];dx)}^2.  \lb{3.55m}
\end{align}
Thus, combining \eqref{3.55i} and \eqref{3.55k}--\eqref{3.55m} one infers
\begin{equation}
\lim_{n\to\infty}\|(H-\lambda I)f_n\|_2 =0, \lb{3.55n}
\end{equation}
and hence $\lambda\in\sigma_{\rm ap}(H)=\sigma(H)$ by \eqref{5.12} and
\eqref{3.41}. 

Relation \eqref{3.45a} is clear from \eqref{3.36} and the fact that by
\eqref{1.3.19} there exists a distributional solution
$\psi((E_m,0),\cdot,x_0)\in L^\infty(\bbR;dx)$ of $H\psi=E_m\psi$ for
all $m=0,\dots,2n$.

Finally, $\sigma(H)$ contains no isolated points since those would
necessarily be essential singularities of the resolvent of $H$, as $H$
has no eigenvalues by \eqref{3.40} (cf.\ \cite[Sect.\ III.6.5]{Ka80}).
An explicit investigation of the Green's function of $H$ reveals at
most a square root singularity at the points $\{E_m\}_{m=0}^{2n}$ and
hence excludes the possibility of an essential singularity of
$(H-zI)^{-1}$.  
\end{proof}

In the special self-adjoint case where $V$ is real-valued, the result
\eqref{3.45} is equivalent to the vanishing of the Lyapunov exponent of
$H$ which characterizes the (purely absolutely continous) spectrum of $H$ 
as discussed by Kotani \cite{Ko84}, \cite{Ko85}, \cite{Ko87a}, \cite{Ko97}
(see also \cite[p.\ 372]{CL90}). In the case where $V$ is periodic and
complex-valued, this has also been studied by Kotani \cite{Ko97}.

The explicit formula for $\Sigma$ in \eqref{3.34} permits a qualitative
description of the spectrum of $H$ as follows. We recall \eqref{3.17}
and write
\begin{equation}
\f{d}{dz}\big\langle g(z,\cdot)^{-1}\big\rangle 
=-2\langle g(z,\cdot)\rangle
=-i \f{\prod_{j=1}^n \big(z-\wti \lambda_j\big)}{\big(\prod_{m=0}^{2n}
(z-E_m)\big)^{1/2}}, \quad z\in\Pi,  \lb{3.90}
\end{equation}
for some constants 
\begin{equation}
\{\wti\lambda_j\}_{j=1}^{n}\subset\bbC. \lb{3.88}
\end{equation}
As in similar situations before, \eqref{3.90} extends to either side of
the cuts in $\calC$ by continuity with respect to $z$.
 
\begin{theorem}  \lb{t3.14} 
Assume Hypothesis \ref{h3.3}. Then the spectrum $\sigma(H)$ of $H$ has
the following properties: \\
$(i)$ $\sigma(H)$ is contained in the semi-strip
\begin{equation}
\sigma(H)\subset \{z\in\bbC\,|\, \Im(z)\in [M_1,M_2], \, \Re(z)\geq
M_3\},  \lb{3.56}
\end{equation}
where
\begin{equation}
M_1=\inf_{x\in\bbR}[\Im(V(x))], \quad 
M_2=\sup_{x\in\bbR}[\Im(V(x))], \quad M_3=
\inf_{x\in\bbR}[\Re(V(x))]. \lb{3.57}
\end{equation}
$(ii)$ $\sigma(H)$ consists of finitely many simple analytic arcs and
one simple semi-infinite arc. These analytic arcs may only end at the
points $\wti\lambda_1,\dots,\wti\lambda_n$, $E_0,\dots,E_{2n}$, and
at infinity. The semi-infinite arc, $\sigma_\infty$, asymptotically
approaches the half-line 
$L_{\langle V\rangle}=\{z\in\bbC \,|\, z=\langle V\rangle +x, \, x\geq
0\}$ in the following sense: asymptotically, 
$\sigma_\infty$ can be parameterized by
\begin{equation}
\sigma_\infty=\big\{z\in\bbC \,\big|\, z=R+i\,\Im(\langle
V\rangle) +\Oh\big(R^{-1/2}\big) 
\text{ as $R\uparrow\infty$}\big\}.  \lb{3.85}
\end{equation} 
$(iii)$ Each $E_m$, $m=0,\dots,2n$, is met by 
at least one of these arcs. More precisely, a particular $E_{m_0}$ is
hit by precisely $2N_0+1$ analytic arcs, where $N_0\in\{0,\dots,n\}$
denotes the number of $\wti\lambda_j$ that coincide with $E_{m_0}$.
Adjacent arcs meet at an angle $2\pi/(2N_0+1)$ at $E_{m_0}$. $($Thus,
generically, $N_0=0$ and precisely one arc hits $E_{m_0}$.$)$ \\
$(iv)$ Crossings of spectral arcs are permitted. This phenomenon and takes
place precisely when for a particular 
$j_0\in\{1,\dots,n\}$, $\wti\lambda_{j_0}\in\sigma(H)$ such that
\begin{equation}
\Re\big(\big\langle g(\wti\lambda_{j_0},\cdot)^{-1}\big\rangle\big)=0 
\, \text{ for some $j_0\in\{1,\dots,n\}$ with $\wti\lambda_{j_0}\notin
\{E_m\}_{m=0}^{2n}$}. \lb{3.85a}
\end{equation}
In this case $2M_0+2$ analytic arcs are converging toward
$\wti\lambda_{j_0}$, where $M_0\in\{1,\dots,n\}$ denotes the number of
$\wti\lambda_j$ that coincide with $\wti\lambda_{j_0}$. Adjacent arcs meet
at an angle $\pi/(M_0+1)$ at $\wti\lambda_{j_0}$. $($Thus, generically,
$M_0=1$ and two arcs cross at a right angle.$)$ \\
$(v)$ The resolvent set $\bbC\backslash\sigma (H)$ of $H$ is
path-connected.
\end{theorem}
\begin{proof} Item $(i)$ follows from \eqref{3.39a} and \eqref{3.41} by
noting that 
\begin{equation}
(f,Hf)=\|f'\|^2+ (f,\Re(V)f)+i(f,\Im(V)f), \quad f\in H^{2,2}(\bbR).
\lb{3.86}
\end{equation}
To prove $(ii)$ we first introduce the meromorphic differential of the
second kind
\begin{equation}
\Omega^{(2)} = \langle g(P,\cdot)\rangle dz= 
\f{i\langle F_n(z,\cdot)\rangle dz}{2y}=\f{i}{2}\f{\prod_{j=1}^n
\big(z-\wti \lambda_j\big) dz}{R_{2n+1}(z)^{1/2}},
\quad  P=(z,y)\in \calK_n\backslash\{\Pinf\} \lb{3.87}
\end{equation}
(cf.\ \eqref{3.88}). Then, by Lemma \ref{l3.4}, 
\begin{equation}
\big\langle g(P,\cdot)^{-1}\big\rangle=
-2\int_{Q_0}^P \Omega^{(2)} + 
\big\langle g(Q_0,\cdot)^{-1}\big\rangle, \quad  
P\in \calK_n\backslash\{\Pinf\} \lb{3.89}
\end{equation}
for some fixed $Q_0\in \calK_n\backslash\{\Pinf\}$, is holomorphic on
$\calK_n\backslash\{\Pinf\}$. By \eqref{3.90}, \eqref{3.88}, the
characterization \eqref{3.45} of the spectrum,
\begin{equation}
\sigma(H) = \big\{\lambda\in\bbC\,\big|\, \Re\big(\big\langle
g(\lambda,\cdot)^{-1}\big\rangle\big)=0\big\}, \lb{3.91}
\end{equation}
and the fact that $\Re\big(\big\langle
g(z,\cdot)^{-1}\big\rangle\big)$ is a harmonic function on the cut
plane $\Pi$, the spectrum $\sigma(H)$ of $H$ consists
of analytic arcs which may only end at the points
$\wti\lambda_1,\dots,\wti\lambda_n$, $E_0,\dots,E_{2n}$, and possibly
tend to infinity. (Since $\sigma(H)$ is independent of the chosen set
of cuts, if a spectral arc crosses or runs along a part of one of the
cuts in $\calC$, one can slightly deform the original set of cuts to
extend an analytic arc along or across such an original cut.) To study the
behavior of spectral arcs near infinity we first note that 
\begin{equation}
g(z,x)\underset{|z|\to\infty}{=}\f{i}{2 z^{1/2}}+\f{i}{4z^{3/2}}V(x)+
\Oh\big(|z|^{-3/2}\big), \lb{3.92}
\end{equation}
combining \eqref{1.2.4}, \eqref{1.2.9}, \eqref{1.2.12}, and
\eqref{3.13}. Thus, one computes
\begin{equation}
g(z,x)^{-1}\underset{|z|\to\infty}{=}-2iz^{1/2} +\f{i}{z^{1/2}}V(x)+
\Oh\big(|z|^{-3/2}\big) \lb{3.93}
\end{equation}
and hence
\begin{equation}
\big\langle g(z,\cdot)^{-1}\big\rangle
\underset{|z|\to\infty}{=}-2iz^{1/2} +\f{i}{z^{1/2}}\langle V\rangle +
\Oh\big(|z|^{-3/2}\big). \lb{3.94}
\end{equation}
Writing $z=R e^{i\varphi}$ this yields
\begin{equation}
0=\Re\big(\big\langle
g(z,\cdot)^{-1}\rangle\big)\underset{R\to\infty}{=} 
2\Im\big\{R^{1/2}e^{i\varphi/2}
-2^{-1}R^{-1/2}e^{-i\varphi/2}\langle V\rangle
+\Oh\big(R^{-3/2}\big)\big\}  \lb{3.95}
\end{equation}
implying
\begin{equation}
\varphi\underset{R\to\infty}{=} \Im(\langle
V\rangle)R^{-1}+\Oh\big(R^{-3/2}\big)  \lb{3.96}
\end{equation}
and hence \eqref{3.85}. In particular, there is precisely one analytic
semi-infinite arc $\sigma_\infty$ that tends to infinity and  
asymptotically approaches the half-line $L_{\langle V\rangle}$. This
proves item $(ii)$.

To prove $(iii)$ one first recalls that by Theorem \ref{t3.12} the spectrum
of $H$ contains no isolated points. On the other hand, since 
$\{E_m\}_{m=0}^{2n}\subset\sigma(H)$ by \eqref{3.45a}, one concludes that
at least one spectral arc meets each $E_m$, $m=0,\dots,2n$. Choosing
$Q_0=(E_{m_0},0)$ in \eqref{3.89} one obtains
\begin{align}
&\big\langle g(z,\cdot)^{-1}\big\rangle =
-2\int_{E_{m_0}}^z dz' \, \langle g(z',\cdot)\rangle + 
\big\langle g(E_{m_0},\cdot)^{-1}\big\rangle  \no \\
&= -i\int_{E_{m_0}}^z dz' \f{\prod_{j=1}^n 
\big(z'-\wti\lambda_j\big)}{\big(\prod_{m=0}^{2n} 
(z'-E_m) \big)^{1/2}}+ 
\big\langle g(E_{m_0},\cdot)^{-1}\big\rangle \no \\ 
&\underset{z\to E_{m_0}}{=} -i\int_{E_{m_0}}^z dz'\,
(z'-E_{m_0})^{N_0-(1/2)}[C+\Oh(z'-E_{m_0})]+ 
\big\langle g(E_{m_0},\cdot)^{-1}\big\rangle \lb{3.97} \\ 
&\underset{z\to E_{m_0}}{=}
-i[N_0+(1/2)]^{-1}(z-E_{m_0})^{N_0+(1/2)}[C+\Oh(z-E_{m_0})]+ 
\big\langle g(E_{m_0},\cdot)^{-1}\big\rangle, \no \\
& \hspace*{11.05cm} z\in \Pi \no
\end{align}
for some $C=|C|e^{i\varphi_0}\in\bbC\backslash\{0\}$. Using 
\begin{equation}
\Re\big(\big\langle g(E_{m},\cdot)^{-1}\big\rangle\big) =0, \quad
m=0,\dots,2n,  \lb{3.98}
\end{equation}
as a consequence of \eqref{3.45a}, $\Re\big(\big\langle
g(z,\cdot)^{-1}\big\rangle\big)=0$ and $z=E_{m_0}+\rho e^{i\varphi}$
imply
\begin{equation}
0\underset{\rho\downarrow 0}{=}
\sin[(N_0+(1/2))\varphi+\varphi_0]\rho^{N_0+(1/2)}[|C|+\Oh(\rho)].
\lb{3.99}
\end{equation}
This proves the assertions made in item $(iii)$.

To prove $(iv)$ it suffices to refer to \eqref{3.90} and to note that
locally, $d \big\langle g(z,\cdot)^{-1}\big\rangle/dz$ behaves like 
$C_0(z-\wti\lambda_{j_0})^{M_0}$ for some $C_0\in\bbC\backslash\{0\}$ in
a sufficiently small neighborhood of $\wti\lambda_{j_0}$.

Finally we will show that all arcs are simple (i.e., do not 
self-intersect each other). Assume that the spectrum of $H$ contains a
simple closed loop $\gamma$, $\gamma\subset\sigma(H)$. Then 
\begin{equation}
\Re\big(\big\langle g(P,\cdot)^{-1}\big\rangle\big)= 0, \quad
P\in\Gamma, \lb{3.100}
\end{equation}
where the closed simple curve $\Gamma\subset\calK_n$ denotes the lift
of $\gamma$ to $\calK_n$, yields the contradiction 
\begin{equation}
\Re\big(\big\langle g(P,\cdot)^{-1}\big\rangle\big)= 0 \, \text{ for 
all $P$ in the interior of $\Gamma$} \lb{3.101}
\end{equation}
by Corollary 8.2.5 in \cite{Be86}. Therefore, since there are no
closed loops in $\sigma(H)$ and precisely one semi-infinite arc tends to
infinity, the resolvent set of $H$ is connected and hence 
path-connected, proving $(v)$.
\end{proof}

\begin{remark} \lb{r4.10}
For simplicity we focused on $L^2(\bbR;dx)$-spectra thus far. However, 
since $V\in L^\infty(\bbR;dx)$, $H$ in $L^2(\bbR;dx)$ is the generator of
a $C_0$-semigroup $T(t)$ in $L^2(\bbR;dx)$, $t>0$, whose integral kernel
$T(t,x,x')$ satisfies the Gaussian upper bound (cf., e.g., \cite{AE97})
\begin{equation}
\big|T(t,x,x')\big|\leq C_1 t^{-1/2}e^{C_2t} e^{-C_3|x-x'|^2/t}, \quad 
t>0, \; x,x'\in\bbR \; \lb{4.93}
\end{equation}
for some $C_1>0$, $C_2\geq 0$, $C_3>0$. Thus, $T(t)$ in $L^2(\bbR;dx)$
defines, for $p\in [1,\infty)$, consistent $C_0$-semigroups
$T_p(t)$ in $L^p(\bbR;dx)$ with generators denoted by $H_p$ (i.e.,
$H=H_2$, $T(t)=T_2(t)$, etc.). Applying Theorem\ 1.1 of Kunstman 
\cite{Ku98} one then infers the $p$-independence of the spectrum, 
\begin{equation}
\sigma(H_p)=\sigma(H), \quad p\in [1,\infty). \lb{4.94}
\end{equation}
Actually, since $\bbC\backslash\sigma (H)$ is connected by 
Theorem \ref{t3.14}\, (v), \eqref{4.94} also follows from Theorem\ 4.2 of
Arendt \cite{Ar94}.
\end{remark}

Of course, these results apply to the special case of algebro-geometric
complex-valued periodic potentials (see \cite{Bi86a}, \cite{Bi86b},
\cite{We98}, \cite{We98a}) and we briefly point out the corresponding
connections between the algebro-geometric approach and standard Floquet
theory in Appendix \ref{sC}. But even in this special case, items (iii) and
(iv) of Theorem \ref{t3.14} provide additional new details on the nature of
the spectrum of $H$.  We briefly illustrate the results of this section in
Example \ref{eC.1} of Appendix \ref{sC}.

The methods of this paper extend to the case of algebro-geometric
non-self-adjoint second order finite difference (Jacobi) operators
associated with the Toda lattice hierarchy and to the case of Dirac-type
operators related to the focusing nonlinear Schr\"odinger hierarchy.
Moreover, they extend to the infinite genus limit
$n\to\infty$ using the approach in
\cite{Ge01}. This will be studied elsewhere.

\appendix
\section{Hyperelliptic curves and their theta functions} \lb{sA}
\renewcommand{\theequation}{A.\arabic{equation}}
\renewcommand{\thetheorem}{A.\arabic{theorem}}
\setcounter{theorem}{0}
\setcounter{equation}{0}

We provide a brief summary of some of the fundamental notations needed
{}from the theory of hyperelliptic Riemann surfaces.  More details can be
found in some of the standard textbooks \cite{FK92} and \cite{Mu84}, as
well as in monographs dedicated to integrable systems such as \cite[Ch.\
2]{BBEIM94}, \cite[App.\ A, B]{GH03}. In particular, the following
material is taken from \cite[App.\ A, B]{GH03}.

Fix $n\in\bbN$. We intend to describe the hyperelliptic 
Riemann surface $\calK_n$ of genus $n$ of the KdV-type curve
(\ref{1.2.17}), associated with the polynomial
\begin{align}
\begin{split}
&\calF_n(z,y)=y^2-R_{2n+1}(z)=0, \lb{b1} \\
&R_{2n+1}(z)=\prod_{m=0}^{2n}(z-E_m),
\quad \{E_m\}_{m=0}^{2n}\subset\bbC. 
\end{split}
\end{align}
To simplify the discussion we will assume that the affine part of
$\calK_n$ is nonsingular, that is, we suppose that 
\begin{equation}
E_m \neq E_{m'} \text{ for } m\neq m', \; m,m'=0,\dots,2n \lb{B1}
\end{equation}
throughout this appendix. Introducing an appropriate set of 
(nonintersecting) cuts $\calC_j$ joining $E_{m(j)}$ and $E_{m^\prime(j)}$,
$j=1,\dots,n$, and $\calC_{n+1}$, joining $E_{2n}$ and $\infty$, we
denote 
\begin{equation}
\calC=\bigcup_{j=1}^{n+1} \calC_j, \quad 
\calC_j\cap\calC_k=\emptyset, \quad j\neq k.\lb{b2}
\end{equation}
Define the cut plane $\Pi$ by 
\begin{equation}
\Pi=\bbC\backslash\calC, \lb{b3}
\end{equation}
and introduce the holomorphic function
\begin{equation}
R_{2n+1}(\cdot)^{1/2}\colon \Pi\to\bbC, \quad 
z\mapsto \bigg(\prod_{m=0}^{2n}(z-E_m) \bigg)^{1/2}\lb{b4}
\end{equation}
on $\Pi$ with an appropriate choice of the square root branch 
in \eqref{b4}. Define
\begin{equation}
\calM_{n}=\{(z,\sigma R_{2n+1}(z)^{1/2}) \mid
z\in\bbC,\; \sigma\in\{1,-1\}
\}\cup \{P_{\infty}\} \label{b5}
\end{equation}
by extending $R_{2n+1}(\cdot)^{1/2}$ to $\calC$. The
hyperelliptic curve $\calK_n$ is then the set
$\calM_{n}$ with its natural complex structure obtained
upon gluing the two sheets of $\calM_{n}$
crosswise along the cuts. The set of branch points
$\calB(\calK_n)$ of $\calK_n$ is given by
\begin{equation}
\calB(\calK_n)=\{(E_m,0)\}_{m=0}^{2n}. \lb{5a}
\end{equation}
Points $P\in\calK_n\backslash\{P_\infty\}$ are denoted by
\begin{equation}
P=(z,\sigma R_{2n+1}(z)^{1/2})=(z,y),  \lb{b15a} 
\end{equation}
where $y(P)$ denotes the meromorphic function
on $\calK_n$ satisfying $\calF_n(z,y)=y^2-R_{2n+1}(z)=0$ and  
\begin{align}
& \quad y(P)\underset{\zeta\to
0}{=}\bigg(1-\f12\bigg(\sum_{m=0}^{2n}E_m\bigg)\zeta^2
+\Oh(\zeta^4)\bigg)\zeta^{-2n-1} \text{ as } P\to P_\infty, \lb{b16} \\ 
& \hspace*{5.9cm} \zeta=\sigma'/z^{1/2}, \, \sigma'\in\{1,-1\} \no
\end{align}
(i.e., we abbreviate $y(P)=\sigma R_{2n+1}(z)^{1/2}$). 
Local coordinates near $P_0=(z_0,y_0)\in\calK_n\backslash
\{\calB(\calK_n)\cup \{P_{\infty}\}\}$ are
given by $\zeta_{P_0}=z-z_0$, near $P_{\infty}$ by
$\zeta_{P_{\infty_\pm}}=1/z^{1/2}$, and near branch points
$(E_{m_0},0)\in\calB(\calK_n)$ by $\zeta_{(E_{m_0},0)}=(z-E_{m_0})^{1/2}$. 
The compact hyperelliptic Riemann surface $\calK_n$ resulting in this
manner has topological genus $n$.  

Moreover, we introduce the holomorphic sheet exchange map (involution)
\begin{equation}
*\colon\calK_n\to\calK_n, \quad P=(z,y)\mapsto P^*=(z,-y), 
\; P_\infty\mapsto P^*_\infty=P_\infty  \lb{b17}
\end{equation}
and the two meromorphic projection maps
\begin{equation}
\tilde\pi\colon\calK_n\to\bbC\cup\{\infty\}, \quad P=(z,y)\mapsto z, \; 
P_\infty\mapsto \infty \lb{b18}
\end{equation}
and
\begin{equation}
y\colon\calK_n\to\bbC\cup\{\infty\}, \quad P=(z,y)\mapsto y, \; 
P_\infty\mapsto \infty.  \lb{b19}
\end{equation}
The map $\tilde\pi$ has a pole of order $2$ at $P_\infty$, and $y$ has a
pole of order $2n+1$ at
$P_\infty$.  Moreover, 
\begin{equation}
\tilde\pi(P^*)=\tilde\pi(P), \quad y(P^*)=-y(P), \quad P\in\calK_n.
\lb{b20}
\end{equation}
Thus $\calK_n$ is a two-sheeted branched covering of the Riemann sphere 
$\bbC\bbP^1$ ($\cong\bbC\cup\{\infty\}$) branched at the $2n+2$
points $\{(E_m,0)\}_{m=0}^{2n}, P_\infty$. 

We introduce the upper and lower sheets $\Pi_{\pm}$ by 
\begin{equation}
\Pi_{\pm}=\{(z,\pm  R_{2n+1}(z)^{1/2})\in \calM_n \mid z\in\Pi\}
\lb{b22}
\end{equation}
and the associated charts
\begin{equation}
\zeta_\pm \colon \Pi_\pm\to \Pi, \quad P\mapsto z.\lb{b23}
\end{equation}

Next, let $\{a_j,b_j\}_{j=1}^n$ be a homology basis for
$\calK_n$ with intersection matrix of the cycles satisfying
\begin{equation}
a_j\circ b_k=\delta_{j,k}, \quad a_j\circ a_k=0, \quad 
b_j\circ b_k=0, \quad j,k=1,\dots,n. \lb{aa16a}
\end{equation}
Associated with the homology basis $\{a_j, b_j\}_{j=1}^n$ we
also recall the canonical dissection of $\calK_n$ along its cycles yielding
the simply connected interior $\hatt \calK_n$ of the
fundamental polygon $\partial {\hatt \calK}_n$ given by
\begin{equation}
\partial  {\hatt \calK}_n =a_1 b_1 a_1^{-1} b_1^{-1}
a_2 b_2 a_2^{-1} b_2^{-1} \cdots a_n^{-1} b_n^{-1}. \lb{aa19}
\end{equation}
Let $\calM (\calK_n)$ and $\calM^1 (\calK_n)$ denote the
set of meromorphic
functions (0-forms) and meromorphic differentials (1-forms)
on $\calK_n$, respectively. The residue of a meromorphic differential
$\nu\in \calM^1 (\calK_n)$ at a
point $Q \in \calK_n$ is defined by
\begin{equation}
\text{res}_{Q}(\nu)
=\frac{1}{2\pi i} \int_{\gamma_{Q}} \nu,
\lb{a33}
\end{equation}
where $\gamma_{Q}$ is a counterclockwise oriented
smooth simple closed
contour encircling $Q$ but no other pole of
$\nu$.  Holomorphic
differentials are also called Abelian differentials
of the first kind. Abelian differentials of the
second kind $\omega^{(2)} \in \calM^1 (\calK_n)$ are characterized
by the property that all their residues vanish.  They will 
usually be normalized by demanding that all their $a$-periods
vanish, that is,
\begin{equation}
\int_{a_j} \omega^{(2)} =0, \quad  j=1,\dots,n.
\lb{a34}
\end{equation}
If $\omega_{P_1, m}^{(2)}$ is a differential of the second kind on
$\calK_n$ whose only pole is $P_1 \in \hatt \calK_n$ with principal part
$\zeta^{-m-2}\,d\zeta$, $m\in\bbN_0$, near $P_1$ and $\omega_j =
\big(\sum_{q=0}^\infty d_{j,q} (P_1) \zeta^q\big) d\zeta$ near $P_1$,
then
\begin{equation}
\frac{1}{2\pi i} \int_{b_j} \omega_{P_1, m}^{(2)} =
 \frac{d_{j,m} (P_1)}{m+1}, \quad m\in\bbN_0, \; j=1,\dots,n. \lb{a35}
\end{equation}

Using the local chart near $P_\infty$, one verifies that $d z/y$ is a
holomorphic differential on $\calK_n$ with zeros of order $2(n-1)$ at 
$P_\infty$ and hence
\begin{equation}
\eta_j=\frac{z^{j-1}d z}{y}, \quad j=1,\dots,n,
\lb{b24}
\end{equation}
form a basis for the space of holomorphic differentials on $\calK_n$.
Upon introduction of the invertible matrix $C$ in $\bbC^n$,
\begin{align}
C & =\big(C_{j,k}\big)_{j,k=1,\dots,n}, \quad C_{j,k}
= \int_{a_k} \eta_j, \lb{b25} \\
\underline{c} (k) & = (c_1(k), \dots,
c_n(k)), \quad c_j (k) =\big(C^{-1}\big)_{j,k}, \quad j,k=1,\dots,n,
\lb{b25a}
\end{align}
the normalized differentials $\ome_j$ for $j=1,\dots,n$, 
\begin{equation}
\ome_j = \sum_{\ell=1}^n c_j (\ell) \eta_\ell,
\quad \int_{a_k} \ome_j =
\delta_{j,k}, \quad j,k=1,\dots,n, \lb{b26}
\end{equation}
form a canonical basis for the space of
holomorphic differentials on $\calK_n$.  

In the chart $(U_{P_\infty}, \zeta_{P_\infty})$ induced by 
$1/\tilde\pi^{1/2}$ near $P_\infty$ one infers,
\begin{align}
{\ul \omega} & = (\omega_1,\dots,\omega_n)= 
-2 \bigg( \sum_{j=1}^n \f{\uc (j)
\zeta^{2(n-j)}}{\big(\prod_{m=0}^{2n}
(1-\zeta^2 E_m) \big)^{1/2}} \bigg) d\zeta \lb{b27} \\
& = -2 \bigg( \uc (n) +\bigg( \frac12 \uc (n)
\sum_{m=0}^{2n} E_m +\uc
(n-1) \bigg) \zeta^2 + \Oh(\zeta^4) \bigg)d\zeta \text{ as $P\to\Pinf$,}
\no \\
& \hspace*{7.1cm} \zeta=\sigma/z^{1/2}, \, \sigma\in\{1,-1\}, \no
\end{align}
where $\ul E=(E_0,\dots, E_{2n})$ and we used \eqref{b16}. 
Given \eqref{b27}, one computes for the vector
$\ul{U}_{0}^{(2)}$ of $b$-periods of $\omega_{P_\infty,0}^{(2)}/(2\pi
i)$, the normalized differential of the second kind, holomorphic on
$\calK_n\backslash\{\Pinf\}$, with principal part
$\zeta^{-2}d\zeta/(2\pi i)$, 
\begin{equation}
\ul{U}_{0}^{(2)}=\big({U}_{0,1}^{(2)},\dots,{U}_{0,n}^{(2)}\big), \quad 
{U}_{0,j}^{(2)}=\f{1}{2\pi i}\int_{b_j} \omega_{P_\infty,0}^{(2)}
=-2 c_j(n), \; j=1,\dots,n. \lb{b27b}  
\end{equation}

Next, define the matrix $\tau=\big(\tau_{j,\ell}\big)_{j,\ell=1}^n$ by
\begin{equation}
\tau_{j,\ell}=\int_{b_j}\omega_\ell, \quad j,\ell=1,
\dots,n. \label{b8}
\end{equation}
Then
\begin{equation}
\Im(\tau)>0, \quad \text{and} \quad \tau_{j,\ell}=\tau_{\ell,j},
\quad j,\ell =1,\dots,n.  \lb{a18a}
\end{equation}
Associated
with $\tau$ one introduces the period lattice
\begin{equation}
L_n = \{ \ul z \in\bbC^n \mid \ul z = \ul m + \ul n\tau,
\; \ul m, \ul n \in\bbZ^n\}
\lb{a28}
\end{equation}
and the Riemann theta function associated with $\calK_n$ and
the given homology basis $\{a_j,b_j\}_{j=1,\dots,n}$,
\begin{equation}
\theta(\ul z)=\sum_{\ul n\in\bbZ^n}\exp\big(2\pi
i(\ul n,\ul z)+\pi
i(\ul n, \ul n\tau)\big),
\quad \ul z\in\bbC^n, \label{b9}
\end{equation}
where $(\ul u, \ul v)= \ol{\ul u}\,\ul v^\top
=\sum_{j=1}^n \overline{u}_j v_j$
denotes the
scalar product
in $\bbC^n$. It has the fundamental properties
\begin{align}
& \theta(z_1, \ldots, z_{j-1}, -z_j, z_{j+1},
\ldots, z_n) =\theta
(\ul z), \lb{a27}\\
& \theta (\ul z +\ul m + \ul n\tau)
=\exp \big(-2 \pi i (\ul n,\ul z) -\pi i (\ul n, 
\ul n\tau) \big) \theta (\ul z), \quad \ul m, \ul n \in\bbZ^n.
\lb{aa51}
\end{align}

Next we briefly study some consequences of a change of homology basis. 
Let 
\begin{equation}
\{a_1,\dots,a_n,b_1,\dots,b_n\} \lb{aa22A}
\end{equation}
be a canonical homology basis on $\calK_n$ with intersection matrix 
satisfying \eqref{aa16a} and let
\begin{equation}
\{a'_1,\dots,a'_n,b'_1,\dots,b'_n\} \lb{aa22B}
\end{equation}  
be a homology basis on $\calK_n$ related to each other by
\begin{equation}
\begin{pmatrix} {\ul a'}^\top \\ {\ul b'}^\top \end{pmatrix}
= X \begin{pmatrix} \ul a^\top \\ \ul b^\top \end{pmatrix}, \lb{aa22C}
\end{equation}
where 
\begin{align} 
\ul a^\top &=(a_1,\dots,a_n)^\top, \;\;\;\;\, \ul b^\top 
=(b_1,\dots,b_n)^\top, \no \\ 
{\ul a'}^\top &=(a'_1,\dots,a'_n)^\top, \quad 
{\ul b'}^\top =(b'_1,\dots,b'_n)^\top, \lb{aa22D} \\
X&=\begin{pmatrix} A & B \\ C & D \end{pmatrix}, \lb{aa22E}
\end{align}
with $A,B,C$, and $D$ being $n\times n$ matrices with integer entries.
Then \eqref{aa22B} is also a canonical homology basis on $\calK_n$ with 
intersection matrix satisfying \eqref{aa16a} if and only if 
\begin{equation}
X \in \Sp(n,\bbZ), \lb{aa22F}
\end{equation}
where
\begin{equation}
\Sp(n,\bbZ)=\left\{X=\begin{pmatrix} A & B \\ C & D
\end{pmatrix}\,\bigg|\, X\begin{pmatrix} 0 & I_n
\\  -I_n & 0 \end{pmatrix}X^\top=\begin{pmatrix} 0 &
I_n \\ -I_n & 0 \end{pmatrix}, \, \det(X)=1\right\} \lb{aa22G}
\end{equation}
denotes the symplectic modular group (here  $A,B,C,D$ in $X$ are again
$n\times n$ matrices with integer entries). If
$\{\omega_j\}_{j=1}^n$ and $\{\omega'_j\}_{j=1}^n$ are the normalized
bases of holomorphic differentials corresponding to the canonical
homology bases \eqref{aa22A} and \eqref{aa22B}, with $\tau$ and $\tau'$
the associated $b$ and $b'$-periods of $\omega_1,\dots,\omega_n$ and
$\omega'_1,\dots,\omega'_n$, respectively, one computes
\begin{equation}
{\ul \omega'}={\ul\omega}(A+B\tau)^{-1}, \quad 
\tau'=(C+D\tau)(A+B\tau)^{-1}, \lb{aa22I}
\end{equation}
where ${\ul\omega}=(\omega_1,\dots,\omega_n)$ and  
${\ul\omega'}=(\omega'_1,\dots,\omega'_n)$. 

Fixing a base point $Q_0\in\calK_n\backslash\{\Pinf\}$, one denotes by
$J(\calK_n) = \bbC^n/L_n$ the Jacobi variety of $\calK_n$,
and defines the Abel map $\underline{A}_{Q_0}$ by
\begin{equation}
\underline{A}_{Q_0} \colon \calK_n \to J(\calK_n), \quad
\underline{A}_{Q_0}(P)=
\bigg(\int_{Q_0}^P \omega_1,\dots,\int_{Q_0}^P \omega_n \bigg)
\pmod{L_n}, \quad P\in\calK_n. \label{aa46}
\end{equation}
Similarly, we introduce
\begin{equation}
\ul \alpha_{Q_0}  \colon
\Div(\calK_n) \to J(\calK_n),\quad
\calD \mapsto \ul \alpha_{Q_0} (\calD)
=\sum_{P \in \calK_n} \calD (P) \ul A_{Q_0} (P),
\label{aa47}
\end{equation}
where $\Div(\calK_n)$ denotes the set of
divisors on $\calK_n$. Here $\calD \colon \calK_n \to \bbZ$
is called a divisor on $\calK_n$ if $\calD(P)\neq0$ for only
finitely many $P\in\calK_n$. (In the main body of this paper
we will choose $Q_0$ to be one of the branch points, i.e.,
$Q_0\in\calB(\calK_n)$, and for simplicity we will always choose
the same path of integration {}from $Q_0$ to $P$ in all Abelian
integrals.) For subsequent use in Remark \ref{raa26a} we also introduce
\begin{align}
\hua_{Q_0} & \colon\hatt{\calK}_n\to\bbC^n, \lb{aa52} \\
&P\mapsto\hua_{Q_0}(P)
=\big(\hatt A_{Q_0,1}(P),\dots,\hatt A_{Q_0,n}(P)\big)
=\bigg(\int_{Q_0}^P\omega_1,\dots,\int_{Q_0}^P\omega_n\bigg) \no 
\end{align}
and
\begin{equation}
\hatt {\ul \al}_{Q_0}  \colon
\Div(\hatt\calK_n) \to \bbC^n, \quad
\calD \mapsto \hatt {\ul \al}_{Q_0} (\calD)
=\sum_{P \in \hatt\calK_n} \calD (P) \hua_{Q_0} (P). \lb{aa52a}
\end{equation}

In connection with divisors on $\calK_n$ we shall employ the
following (additive) notation,
\begin{align} 
&\calD_{Q_0\ul Q}=\calD_{Q_0}+\calD_{\ul Q}, \quad \calD_{\ul
Q}=\calD_{Q_1}+\cdots +\calD_{Q_m}, \lb{A.17} \\
& {\ul Q}=\{Q_1, \dots ,Q_m\} \in \sym^m \calK_n,
\quad Q_0\in\calK_n, \; m\in\bbN, \no
\end{align}
where for any $Q\in\calK_n$,
\begin{equation} \lb{A.18}
\calD_Q \colon  \calK_n \to\bbN_0, \quad
P \mapsto  \calD_Q (P)=
\begin{cases} 1 & \text{for $P=Q$},\\
0 & \text{for $P\in \calK_n\backslash \{Q\}$}, \end{cases}
\end{equation}
and $\sym^m \calK_n$ denotes the $m$th symmetric product of
$\calK_n$. In particular, $\sym^m \calK_n$ can be
identified with
the set of nonnegative
divisors $0 \leq \calD \in \Div(\calK_n)$ of degree $m\in\bbN$.

For $f\in \calM (\calK_n) \backslash \{0\}$ and 
$\omega \in \calM^1 (\calK_n) \backslash \{0\}$ the
divisors of $f$ and $\omega$ are denoted
by $(f)$ and
$(\omega)$, respectively.  Two
divisors $\calD$, $\calE\in \Div(\calK_n)$ are
called equivalent, denoted by
$\calD \sim \calE$, if and only if $\calD -\calE
=(f)$ for some
$f\in\calM (\calK_n) \backslash \{0\}$.  The divisor class
$[\calD]$ of $\calD$ is
then given by $[\calD]
=\{\calE \in \Div(\calK_n)\mid\calE \sim \calD\}$.  We
recall that
\begin{equation}
\deg ((f))=0,\, \deg ((\omega)) =2(n-1),\,
f\in\calM (\calK_n) \backslash
\{0\},\,  \omega\in \calM^1 (\calK_n) \backslash \{0\},
\lb{a38}
\end{equation}
where the degree $\deg (\calD)$ of $\calD$ is given
by $\deg (\calD)
=\sum_{P\in \calK_n} \calD (P)$.  It is customary to call
$(f)$ (respectively,
$(\omega)$) a principal (respectively, canonical)
divisor.

Introducing the complex linear spaces
\begin{align}
\calL (\calD) & =\{f\in \calM (\calK_n)\mid f=0
 \text{ or } (f) \geq \calD\}, \;
r(\calD) =\dim_\bbC \calL (\calD),
\lb{a39}\\
\calL^1 (\calD) & =
 \{ \omega\in \calM^1 (\calK_n)\mid \omega=0
 \text{ or } (\omega) \geq
\calD\},\; i(\calD) =\dim_\bbC \calL^1 (\calD)  \lb{a40}
\end{align}
(with $i(\calD)$ the index of specialty of $\calD$), one
infers that $\deg
(\calD)$, $r(\calD)$, and $i(\calD)$ only depend on
the divisor class
$[\calD]$ of $\calD$.  Moreover, we recall the
following fundamental
facts.

\begin{theorem} \lb{thm1}
Let $\calD \in \Div(\calK_n)$,
$\omega \in \calM^1 (\calK_n) \backslash \{0\}$. Then
\begin{equation}
 i(\calD) =r(\calD-(\omega)), \quad n\in\bbN_0.
\lb{a41}
\end{equation}
The Riemann-Roch theorem reads
\begin{equation}
r(-\calD) =\deg (\calD) + i (\calD) -n+1,
\quad n\in\bbN_0.
\lb{a42}
\end{equation}
By Abel's theorem, $\calD\in \Div(\calK_n)$,
$n\in\bbN$, is principal
if and only if
\begin{equation}
\deg (\calD) =0 \text{ and } \ul \alpha_{Q_0} (\calD)
=\ul{0}.
\lb{a43}
\end{equation}
Finally, assume
$n\in\bbN$. Then $\ul \alpha_{Q_0}
: \Div(\calK_n) \to J(\calK_n)$ is surjective
$($Jacobi's inversion theorem$)$.
\end{theorem}

\begin{theorem} \lb{thm3}
Let $\calD_{\ul Q} \in \sym^n \calK_n$,
$\ul Q=\{Q_1, \ldots, Q_n\}$.  Then
\begin{equation}
1 \leq i (\calD_{\ul Q} ) =s
\lb{a46}
\end{equation}
if and only if there are $s$ pairs of the type $\{P, P^*\}\subseteq 
\{Q_1,\ldots, Q_n\}$ $($this includes, of course, branch points for which
$P=P^*$$)$. Obviously, one has $s\leq n/2$.
\end{theorem}

Next, denote by $\ul \Xi_{Q_0}=(\Xi_{Q_{0,1}}, \dots,
\Xi_{Q_{0,n}})$ the vector of Riemann constants,
\begin{equation}
\Xi_{Q_{0,j}}=\frac12(1+\tau_{j,j})-
\sum_{\substack{\ell=1 \\ \ell\neq j}}^n\int_{a_\ell}
\omega_\ell(P)\int_{Q_0}^P\omega_j,
\quad j=1,\dots,n. \lb{aa55}
\end{equation}

\begin{theorem} \lb{taa17a}
Let $\ul Q =\{Q_1,\dots,Q_n\}\in \sym^n \calK_n$ and
assume $\calD_{\ul Q}$ to be nonspecial, that is,
$i(\calD_{\ul Q})=0$. Then
\begin{equation}
\theta(\ul {\Xi}_{Q_0} -\ul {A}_{Q_0}(P) + \alpha_{Q_0}
(\calD_{\ul Q}))=0 \text{ if and only if }
P\in\{Q_1,\dots,Q_n\}. \lb{aa55a}
\end{equation}
\end{theorem}

\begin{remark} \lb{raa26a}
In Section \ref{s2} we dealt with theta function expressions of the type 
\begin{equation}
\psi(P)=\f{\theta(\ul\Xi_{Q_0}-\ul A_{Q_0}(P)+\ul\alpha_{Q_0}(\calD_1))}
{\theta(\ul\Xi_{Q_0}-\ul A_{Q_0}(P)+\ul\alpha_{Q_0}(\calD_2))}
\exp\bigg(-c \int_{Q_0}^P \Omega^{(2)}\bigg), \quad P\in\calK_n,
\lb{aa76b}
\end{equation}
where $\calD_j\in\sym^n\calK_n$, $j=1,2$, are nonspecial positive
divisors of degree $n$, $c\in\bbC$ is a constant, and $\Omega^{(2)}$ is a
normalized differential of the second kind with a prescribed singularity
at $\Pinf$. Even though we agree to always choose identical paths of
integration {}from $P_0$ to $P$ in all Abelian integrals 
\eqref{aa76b}, this is not sufficient to render $\psi$
single-valued on $\calK_n$. To achieve single-valuedness one needs to
replace $\calK_n$ by its simply connected canonical dissection
$\hatt\calK_n$ and then replace $\ul A_{Q_0}$ and $\ul \alpha_{Q_0}$ 
in \eqref{aa76b} with ${\hua}_{Q_0}$ and $\hatt{\ul
\alpha}_{Q_0}$ as introduced in \eqref{aa52} and \eqref{aa52a}. In 
particular, one regards $a_j,b_j$, $j=1,\dots,n$, as curves (being a 
part of $\partial\hatt\calK_n$, cf. \eqref{aa19}) and not as homology
classes. Similarly, one then replaces $\uxi_{Q_0}$ by \,$\hatt \uxi_{Q_0}$
(replacing $\ul A_{Q_0}$ by ${\hua}_{Q_0}$ in \eqref{aa55}, etc.).  
Moreover, in connection with $\psi$, one introduces the vector of
$b$-periods $\ul U^{(2)}$ of $\Omega^{(2)}$ by
\begin{equation}
\ul U^{(2)}=(U_1^{(2)},\dots,U_n^{(2)}), \quad U_j^{(2)}=\f{1}{2\pi
i}\int_{b_j} \Omega^{(2)}, \quad j=1,\dots,n, \lb{aa76e}
\end{equation}
and then renders $\psi$ single-valued on $\hatt\calK_n$ by requiring 
\begin{equation}
\hatt{\ul\alpha}_{Q_0}(\calD_1)-\hatt{\ul\alpha}_{Q_0}(\calD_2)
=c \,\ul U^{(2)} \lb{aa76f}
\end{equation}
$($as opposed to merely
$\ul\alpha_{Q_0}(\calD_1)-\ul\alpha_{Q_0}(\calD_2)=c \,\ul U^{(2)}
\pmod {L_n}$$)$. Actually, by \eqref{aa51},
\begin{equation}
\hatt{\ul\alpha}_{Q_0}(\calD_1)-\hatt{\ul\alpha}_{Q_0}(\calD_2)
- c\, \ul U^{(2)}\in\bbZ^n, \lb{aa76h}
\end{equation}
suffices to guarantee single-valuedness of $\psi$ on
$\hatt\calK_n$. Without the replacement of $\ul A_{Q_0}$ and  
$\ul \alpha_{Q_0}$ by ${\hua}_{Q_0}$ and $\hatt{\ul \alpha}_{Q_0}$  in
\eqref{aa76b} and without the assumption  
\eqref{aa76f} $($or \eqref{aa76h}$)$, $\psi$ is a
multiplicative $($multi-valued$)$ function on $\calK_n$, and then most
effectively discussed by introducing the notion of characters on
$\calK_n$ $($cf.\ \cite[Sect.\ III.9]{FK92}$)$. For simplicity, we
decided to avoid the latter possibility and throughout this paper will
always tacitly assume \eqref{aa76f} or \eqref{aa76h}.
\end{remark}

\section{Restrictions on $\ul B = i \ul U_0^{(2)}$} \lb{sB}
\renewcommand{\theequation}{B.\arabic{equation}}
\renewcommand{\thetheorem}{B.\arabic{theorem}}
\setcounter{theorem}{0}
\setcounter{equation}{0}

The purpose of this appendix is to prove the result \eqref{2.70}, 
$\ul B=i \ul U^{(2)}_0 \in \bbR^n$, for some choice of homology basis 
$\{a_j, b_j\}_{j=1}^n$ on $\calK_n$ as recorded in Remark \ref{r2.8}.

To this end we first recall a few notions in connection with
periodic meromorphic functions of $p$ complex variables. 

\begin{definition} \lb{dB.1}
Let $p\in\bbN$ and $F\colon\bbC^p\to\bbC\cup\{\infty\}$ be meromorphic
(i.e., a ratio of two entire functions of $p$ complex variables). 
Then, \\
(i) $\ul \omega=(\omega_1,\dots,\omega_p)\in\bbC^p\backslash\{0\}$ is
called a {\it period} of $F$ if
\begin{equation}
F(\ul z+\ul \omega)=F(\ul z) \lb{B.1}
\end{equation}
for all $\ul z\in\bbC^p$ for which $F$ is analytic. The set of all
periods of $F$ is denoted by $\calP_F$. \\
(ii) $F$ is called {\it degenerate} if it depends on less than $p$
complex variables; otherwise, $F$ is called {\it nondegenerate}. 
\end{definition}

\begin{theorem} \lb{tB.2}
Let $p\in\bbN$, $F\colon\bbC^p\to\bbC\cup\{\infty\}$ be meromorphic, and 
$\calP_F$ be the set of all periods of $F$. Then either \\
$(i)$ $\calP_F$ has a finite limit point, \\
or \\
$(ii)$ $\calP_F$ has no finite limit point. \\
In case $(i)$, $\calP_F$ contains {\it infinitesimal periods}
$($i.e., sequences of nonzero periods converging to zero$)$. In
addition, in case $(i)$ each period is a limit point of periods
and hence $\calP_F$ is a perfect set. \\
Moreover, $F$ is degenerate if and only if $F$ admits infinitesimal
periods. In particular, for nondegenerate functions $F$ only
alternative $(ii)$ applies.  
\end{theorem}

Next, let $\ul\omega_q\in\bbC^p\backslash\{0\}$, $q=1,\dots,r$ for
some $r\in\bbN$. Then $\ul\omega_1,\dots,\ul\omega_r$ are called
{\it linearly independent over $\bbZ$ $($resp.\ $\bbR$$)$} if 
\begin{align}
&\nu_1\ul\omega_1+\cdots+\nu_r\ul\omega_r=0, \quad \nu_q\in\bbZ 
\text{ (resp., $\nu_q\in\bbR$)}, \; q=1,\dots,r, \no \\
&\text{implies } \nu_1=\cdots=\nu_r=0. \lb{B.2} 
\end{align}
Clearly, the maximal number of vectors in $\bbC^p$ linearly
independent over $\bbR$ equals $2p$. 

\begin{theorem} \lb{tB.3}
Let $p\in\bbN$. \\
$(i)$ If $F\colon\bbC^p\to\bbC\cup\{\infty\}$ is a nondegenerate
meromorphic function with periods
$\ul\omega_q\in\bbC^p\backslash\{0\}$, $q=1,\dots,r$, $r\in\bbN$, 
linearly independent over $\bbZ$, then $\ul\omega_1,\dots,\ul\omega_r$
are also linearly independent over $\bbR$. In particular, $r\leq 2p$.
\\
$(ii)$ A nondegenerate entire function $F\colon\bbC^p\to\bbC$ cannot
have more than $p$ periods linearly independent over $\bbZ$ $($or
$\bbR$$)$.  
\end{theorem}

For $p=1$, $\exp(z)$, $\sin(z)$ are examples of entire functions with
precisely one period. Any non-constant doubly periodic meromorphic
function of one complex variable is elliptic (and hence has indeed
poles). 

\begin{definition} \lb{dB.4}
Let $p, r\in\bbN$. A system of periods
$\ul\omega_q\in\bbC^p\backslash\{0\}$, $q=1,\dots,r$ of a
nondegenerate meromorphic function
$F\colon\bbC^p\to\bbC\cup\{\infty\}$, linearly independent over $\bbZ$,
is called {\it fundamental} or a {\it basis} of periods for $F$ if every
period $\ul\omega$ of $F$ is of the form
\begin{equation}
\ul\omega =m_1\ul\omega_1+\cdots+m_r\ul\omega_r \, 
\text{ for some $m_q\in\bbZ$, $q=1,\dots,r$.} \lb{B.3}
\end{equation}
\end{definition}

The representation of $\ul\omega$ in \eqref{B.3} is unique since by
hypothesis $\ul\omega_1,\dots,\ul\omega_r$ are linearly independent
over $\bbZ$. In addition, $\calP_F$ is countable in this case. (This
rules out case $(i)$ in Theorem \ref{tB.2} since a perfect set is
uncountable. Hence, one does not have to assume that $F$ is
nondegenerate in Definition \ref{dB.4}.)

This material is standard and can be found, for instance, in
\cite[Ch.\ 2]{Ma92}. 

\vspace*{2mm}

Next, returning to the Riemann theta function $\theta(\ul\cdot)$ in
\eqref{b9}, we introduce the vectors $\{\ul e_j\}_{j=1}^n,
\{\ul\tau_j\}_{j=1}^n
\subset\bbC^n\backslash\{0\}$ by
\begin{equation}
\ul e_j = (0,\dots,0,\underbrace{1}_{j},0,\dots,0), \quad 
\ul \tau_j = \ul e_j \tau, \quad j=1,\dots,n. \lb{B.4}
\end{equation}
Then 
\begin{equation}
\{\ul e_j\}_{j=1}^n \lb{B.5}
\end{equation}
is a basis of periods for the entire (nondegenerate) function 
$\theta(\ul\cdot)\colon\bbC^n\to\bbC$. Moreover, fixing
$k,k'\in\{1,\dots,n\}$, then 
\begin{equation}
\{\ul e_j, \ul\tau_j\}_{j=1}^n \lb{B.6}
\end{equation}
is a basis of periods for the meromorphic function $\partial^2_{z_k
z_{k'}}\ln\big(\theta(\ul\cdot)\big)\colon\bbC^n\to\bbC\cup\{\infty\}$
(cf.\ \eqref{aa51} and \cite[p.\ 91]{FK92}). 

Next, let $\ul A\in\bbC^n$, $\ul D=(D_1,\dots,D_n)\in\bbR^n$, 
$D_j\in\bbR\backslash\{0\}$, $j=1,\dots,n$ and consider
\begin{align}
\begin{split}
f_{k,k'}\colon \bbR\to\bbC, \quad 
f_{k,k'}(x)&=\partial^2_{z_k z_{k'}} 
\ln\big(\theta(\ul A+\ul z)\big)\big|_{\ul z=\ul D x} \lb{B.7} \\
&=\partial^2_{z_k z_{k'}} 
\ln\big(\theta(\ul A+\ul z\diag(\ul D))\big)\big|_{\ul z=(x,\dots,x)}. 
\end{split}
\end{align}
Here $\diag(\ul D)$ denotes the diagonal matrix 
\begin{equation}
\diag(\ul D)= \big(D_j\delta_{j,j'}\big)_{j,j'=1}^n. \lb{B.8}
\end{equation}
Then the quasi-periods $D_j^{-1}$, $j=1,\dots,n$, of $f_{k, k'}$ are in
a one-to-one correspondence with the periods of 
\begin{equation}
F_{k,k'}\colon\bbC^n\to\bbC\cup\{\infty\}, \quad 
F_{k,k'}(\ul z)=\partial^2_{z_k z_{k'}} 
\ln\big(\theta(\ul A+\ul z\diag(\ul D)\big) \lb{B.9}
\end{equation}
of the special type
\begin{equation}
\ul e_j \big(\diag (\ul D)\big)^{-1} =
\big(0,\dots,0,\underbrace{D_j^{-1}}_{j},0,\dots,0\big). \lb{B.10}
\end{equation} 
Moreover,
\begin{equation}
f_{k,k'}(x)=F_{k,k'}(\ul z)|_{\ul z=(x,\dots,x)}, \quad x\in\bbR. 
\lb{B.11}
\end{equation}

\begin{theorem} \lb{tB.5}
Suppose $V$ in \eqref{1.3.61a} $($or \eqref{1.3.IM}$)$ to be
quasi-periodic. Then there exists a homology basis 
$\{\ti a_j, \ti b_j\}_{j=1}^n$ on $\calK_n$ such that the vector
$\wti{\ul B}=i\wti{\ul U}^{(2)}_0$ with $\wti{\ul U}^{(2)}_0$ the
vector of $\ti b$-periods of the corresponding normalized
differential of the second kind, $\wti \omega^{(2)}_{\Pinf,0}$,
satisfies the constraint
\begin{equation}
\wti {\ul B}=i \wti{\ul U}^{(2)}_0 \in \bbR^n. \lb{B.12}
\end{equation}
\end{theorem}
\begin{proof}
By \eqref{b27b}, the vector of $b$-periods $\ul U^{(2)}_0$ associated
with a given homology basis $\{a_j, b_j\}_{j=1}^n$ on $\calK_n$ and the
normalized differential of the 2nd kind, $\omega^{(2)}_{\Pinf,0}$, is
continuous with respect to $E_0,\dots,E_{2n}$. Hence, we may assume in
the following that 
\begin{equation}
B_j\neq 0, \; j=1,\dots,n, \quad \ul B=(B_1,\dots,B_n) \lb{B.13}
\end{equation}
by slightly altering $E_0,\dots,E_{2n}$, if necessary. By comparison
with the Its--Matveev formula \eqref{1.3.IM}, we may write
\begin{align}
\begin{split}
V(x)& =\Lambda_0 -2\partial_x^2 \ln(\theta(\ul A +\ul B x)) \\
&=\Lambda_0 +2\sum_{j,k=1}^n U^{(2)}_{0,j} U^{(2)}_{0,k} 
\partial^2_{z_k z_{j}} 
\ln\big(\theta(\ul A+\ul z)\big)\big|_{\ul z=\ul B x}. \lb{B.14}
\end{split}
\end{align}
Introducing the meromorphic (nondegenerate) function 
$\calV\colon\bbC^n\to\bbC\cup\{\infty\}$ by
\begin{equation}
\calV(\ul z)=\Lambda_0 +2\sum_{j,k=1}^n U^{(2)}_{0,j} U^{(2)}_{0,k} 
\partial^2_{z_k z_{j}} 
\ln\big(\theta(\ul A+\ul z\diag(\ul B))\big), \lb{B.15}
\end{equation}
one observes that
\begin{equation}
V(x)=\calV(\ul z)|_{\ul z=(x,\dots,x)}. \lb{B.16}
\end{equation}
In addition, $\calV$ has a basis of periods 
\begin{equation}
\Big\{\ul e_j \big(\diag(\ul B)\big)^{-1}, 
\ul\tau_j \big(\diag(\ul B)\big)^{-1}\Big\}_{j=1}^n \lb{B.17}
\end{equation}
by \eqref{B.6}, where
\begin{align}
\ul e_j \big(\diag(\ul B)\big)^{-1} &=
\big(0,\dots,0,\underbrace{B_j^{-1}}_{j},0,\dots,0\big), 
\quad j=1,\dots,n, \lb{B.18} \\ 
\ul\tau_j\big(\diag(\ul B)\big)^{-1} &=\big(\tau_{j,1}B_1^{-1},\dots,
\tau_{j,n}B_n^{-1}\big), \quad j=1,\dots,n. \lb{B.19}
\end{align}
By hypothesis, $V$ in \eqref{B.14} is quasi-periodic and hence has 
$n$ real (scalar) quasi-periods. The latter are not necessarily
linearly independent over $\bbQ$ from the outset, but by slightly changing
the locations of branchpoints $\{E_m\}_{m=0}^{2n}$ into, say,
$\{\wti E_m\}_{m=0}^{2n}$, one can assume they are. In particular, since
the period vectors in \eqref{B.17} are linearly independent and the
(scalar) quasi-periods of $V$ are in a one-one correspondence with vector
periods of $\calV$ of the special form \eqref{B.18} (cf.\ \eqref{B.9},
\eqref{B.10}), there exists a homology basis $\{\ti a_j, \ti
b_j\}_{j=1}^n$ on $\calK_n$ such that the vector $\ul {\wti B}=i \wti{\ul
U}^{(2)}_0$ corresponding to the normalized differential of the second
kind, $\wti  \omega^{(2)}_{\Pinf,0}$ and this particular homology basis, is
real-valued. By continuity of $\wti{\ul U}^{2}_0$ with respect to 
$\wti E_0,\dots,\wti E_{2m}$, this proves \eqref{B.12}.
\end{proof}

\begin{remark} \lb{rB.6}
Given the existence of a homology basis with associated real vector
$\wti{\ul B}=i \wti{\ul U}^{(2)}_0$, one can follow the proof of
Theorem\ 10.3.1 in \cite{Le87} and show that each $\mu_j$,
$j=1,\dots,n$, is quasi-periodic with the same quasi-periods as $V$.
\end{remark}

\section{Floquet theory and an explicit example} \lb{sC}
\renewcommand{\theequation}{C.\arabic{equation}}
\renewcommand{\thetheorem}{C.\arabic{theorem}}
\setcounter{theorem}{0}
\setcounter{equation}{0}

In this appendix we discuss the special case of algebro-geometric
complex-valued periodic potentials and we briefly point out the
connections between the algebro-geometric approach and standard Floquet
theory. We then conclude with the explicit genus $n=1$ example which
illustrates both, the algebro-geometric as well as the periodic case. 

We start with the periodic case. Suppose $V$ satisfies
\begin{equation}
V\in CP(\bbR) \, \text{ and for all $x\in\bbR$, } \, V(x+\Omega)=V(x)
\end{equation}
for some period $\Omega>0$. In addition, we suppose that $V$ satisfies
Hypothesis \ref{h3.3}.

Under these assumptions the Riemann surface associated with $V$, which by
Floquet theoretic arguments, in general, would be a two-sheeted Riemann
surface of infinite genus, can be reduced to the compact hyperelliptic
Riemann surface corresponding to $\calK_n$ induced by $y^2=R_{2n+1}(z)$.
Moreover, the corresponding Schr\"odinger operator $H$ is then defined as
in \eqref{3.35} and one introduces the fundamental system of
distributional solutions $c(z,\cdot,x_0)$ and $s(z,\cdot,x_0)$ of
$H\psi=z\psi$ satisfying
\begin{align}
&c(z,x_0,x_0)=s_x(z,x_0,x_0)=1, \\
&c_x(z,x_0,x_0)=s(z,x_0,x_0)=0, \quad z\in\bbC
\end{align}
with $x_0\in\bbR$ a fixed reference point. For each $x, x_0\in\bbR$, 
$c(z,x,x_0)$ and $s(z,x,x_0)$ are entire with respect to $z$. The
monodromy matrix
$\calM(z,x_0)$ is then given by
\begin{equation}
\calM(z,x_0)=\begin{pmatrix} c(z,x_0+\Omega,x_0) & s(z,x_0+\Omega,x_0) \\
c_x(z,x_0+\Omega,x_0) & s_x(z,x_0+\Omega,x_0) \end{pmatrix}, \quad
z\in\bbC
\end{equation}
and its eigenvalues $\rho_\pm(z)$, the Floquet multipliers (which are
$x_0$-independent), satisfy
\begin{equation}
\rho_+(z)\rho_-(z)=1
\end{equation}
since $\det(\calM(z,x_0))=1$. The Floquet discriminant $\Delta(\cdot)$ is
then defined by
\begin{equation}
\Delta(z)=\tr(\calM(z,x_0))/2=
[c(z,x_0+\Omega,x_0)+s_x(z,x_0+\Omega,x_0)]/2
\end{equation} 
and one obtains
\begin{equation}
\rho_\pm (z)=\Delta(z)\pm [\Delta(z)^2-1]^{1/2}.
\end{equation}
We also note that 
\begin{equation}
|\rho_\pm(z)|=1 \, \text{ if and only if } \, \Delta(z)\in [-1,1]. 
\end{equation}
The Floquet solutions $\psi_\pm(z,x,x_0)$, the analog of the functions in
\eqref{3.47}, are then given by 
\begin{align}
&\psi_\pm(z,x,x_0)=c(z,x,x_0)+s(z,x,x_0)[\rho_\pm(z)-c(z,x_0+\Omega,x_0)]
s(z,x_0+\Omega,x_0)^{-1}, \no \\
& \hspace*{7.2cm} z\in\Pi\backslash\{\mu_j(x_0)\}_{j=1,\dots,n}
\end{align}
and one verifies (for $x, x_0 \in\bbR$),
\begin{align}
&\psi_\pm(z,x+\Omega,x_0)=\rho_\pm(z) \psi_\pm(z,x,x_0), \quad 
z\in\Pi\backslash\{\mu_j(x_0)\}_{j=1,\dots,n}, \\
&\psi_+(z,x,x_0)\psi_-(z,x,x_0)=\f{s(z,x+\Omega,x)}{s(z,x_0+\Omega,x_0)},
\quad z\in\bbC\backslash\{\mu_j(x_0)\}_{j=1,\dots,n}, \lb{3.118} \\
& W(\psi_+(z,\cdot,x_0),\psi_-(z,\cdot,x_0))= - 
\f{2[\Delta(z)^2-1]^{1/2}}{s(z,x_0+\Omega,x_0)}, \quad 
z\in\Pi\backslash\{\mu_j(x_0)\}_{j=1,\dots,n},  \lb{3.119} \\ 
&g(z,x)=-\f{s(z,x+\Omega,x)}{2[\Delta(z)^2-1]^{1/2}}
=\f{iF_n(z,x)}{2R_{2n+1}(z)^{1/2}}, \quad z\in\Pi. \lb{3.120}
\end{align}
Moreover, one computes
\begin{align}
\f{d\Delta(z)}{dz}&=-s(z,x_0+\Omega,x_0)\f{1}{2} \int_{x_0}^{x_0+\Omega}
dx \, \psi_+(z,x,x_0)\psi_-(z,x,x_0) \no \\
&=\Omega [\Delta(z)^2-1]^{1/2} \langle g(z,\cdot)\rangle, \quad
z\in\bbC 
\lb{3.121}
\end{align}
and hence
\begin{equation}
\f{d \Delta(z)/dz}{[\Delta(z)^2-1]^{1/2}}
=\f{d}{dz}\big\{\ln\big[\Delta(z)+[\Delta(z)^2-1]^{1/2}\big]\big\}=\Omega 
\langle g(z,\cdot)\rangle, \quad z\in\Pi. \lb{3.122}
\end{equation}
Here the mean value $\langle f\rangle$ of a periodic function $f\in
CP(\bbR)$ of period $\Omega>0$ is simply given by 
\begin{equation}
\langle f\rangle=\f{1}{\Omega} \int_{x_0}^{x_0+\Omega} dx\, f(x),
\lb{3.123}
\end{equation}
independent of the choice of $x_0\in\bbR$. Thus, applying \eqref{3.17} one
obtains 
\begin{align}
& \int_{z_0}^z \f{dz'\, [d\Delta(z')/dz']}{[\Delta(z')^2-1]^{1/2}}=
\ln\bigg(\f{\Delta(z)+[\Delta(z)^2-1]^{1/2}}{\Delta(z_0)
+[\Delta(z_0)^2-1]^{1/2}}\bigg) \no \\
&\quad =\Omega \int_{z_0}^z dz'\, \langle g(z',\cdot)\rangle = 
-(\Omega/2)\big[\big\langle g(z,\cdot)^{-1}\big\rangle 
-\big\langle g(z_0,\cdot)^{-1}\big\rangle\big], \quad z,z_0 \in \Pi
\lb{3.124}
\end{align}
and hence 
\begin{equation}
\ln\big[\Delta(z)+[\Delta(z)^2-1]^{1/2}\big]=-(\Omega/2)\big\langle
g(z,\cdot)^{-1}\big\rangle + C. \lb{3.125}
\end{equation}
Letting $|z|\to \infty$ one verifies that $C=0$ and thus
\begin{equation}
\ln\big[\Delta(z)+[\Delta(z)^2-1]^{1/2}\big]=-(\Omega/2)\big\langle
g(z,\cdot)^{-1}\big\rangle, \quad z\in\Pi. \lb{3.126}
\end{equation}
We note that by continuity with respect to $z$, equations
\eqref{3.119}, \eqref{3.120}, \eqref{3.122}, \eqref{3.124},
and \eqref{3.126} all extend to either side of the set of cuts in
$\calC$.  Consequently,
\begin{equation}
\Delta(z)\in [-1,1] \, \text{ if and only if } \, \Re\big(\big\langle
g(z,\cdot)^{-1}\big\rangle\big)=0. \lb{3.127}
\end{equation}
In particular, our characterization of the spectrum of $H$ in \eqref{3.45}
is thus equivalent to the standard Floquet theoretic characterization of
$H$ in terms of the Floquet discriminant,
\begin{equation}
\sigma(H)= \{\lambda\in\bbC\,|\, \Delta(\lambda)\in [-1,1]\}. \lb{3.128}
\end{equation}
The result \eqref{3.128} was originally proven in \cite{Ro63} and
\cite{Se60} for complex-valued periodic (not necessarily
algebro-geometric) potentials (cf.\ also \cite{Tk64}, and more
recently, \cite{Tk94}, \cite{Tk96}).

We will end this appendix by providing an explicit example of the
simple yet nontrivial genus $n=1$ case which illustrates the periodic
case as well as some of the general results of Sections 
\ref{s2}--\ref{s4} and Appendix \ref{sB}. For more general elliptic
examples we refer to \cite{GW96}, \cite{GW98} and the references
therein.

By $\wp(\cdot)=\wp(\cdot\,|\,\Omega_1,\Omega_3)$ we denote the Weierstrass
$\wp$-function with fundamental half-periods $\Omega_j$, $j=1, 3$,
$\Omega_1>0$, $\Omega_3\in\bbC\backslash\{0\}$, $\Im(\Omega_3) > 0$,
$\Omega_2=\Omega_1+\Omega_3$, and invariants
$g_2$ and $g_3$  (cf.\ \cite[Ch.\ 18]{AS72}). By 
$\zeta(\cdot)=\zeta(\cdot\,| \Omega_1, \Omega_3)$ and
$\sigma(\cdot)=\sigma(\cdot\,| \Omega_1, \Omega_3)$ we denote the
Weierstrass zeta and sigma functions, respectively. We also denote
$\tau=\Omega_3/\Omega_1$ and hence stress that $\Im(\tau)>0$.

\begin{example} \lb{eC.1} 
Consider the genus one ($n=1$) Lam\'e potential
\begin{align}
V(x)&=2\wp(x+\Omega_3) \lb{1.3.76}  \\
&= - 2 \bigg\{\ln\bigg[\theta\bigg(\f{1}{2}+\f{x}{2\Omega_1} 
\bigg)\bigg]\bigg\}''-2\f{\zeta(\Omega_1)}{\Omega_1},
\quad x\in\bbR, \lb{1.3.76a}
\end{align} 
where
\begin{equation}
\theta(z)=\sum_{n\in\bbZ} \exp\big(2\pi inz+\pi in^2\tau\big),
\quad  z\in\bbC, \; \tau=\Omega_3/\Omega_1, \lb{1.3.76b}
\end{equation}
and introduce 
\begin{equation}
L  =  -\f{d^2}{dx^2}+2 \wp(x+\Omega_3), \quad
P_3=-\f{d^3}{dx^3}+3\wp(x+\Omega_3)\f{d}{dx}+\f{3}{2}\wp'(x+\Omega_3). 
\lb{1.3.77} 
\end{equation}
Then one obtains
\begin{equation}
[L,P_3]=0 \lb{1.3.77a}
\end{equation}
which yields the elliptic curve
\begin{align}
&\calK_1\colon\calF_{1}(z,y)=y^2-R_3(z)=y^2-\big(z^3-(g_2/4) z
+(g_3/4)\big)=0, \no  \\
& R_3(z)=\prod_{m=0}^2 (z-E_m)=z^3-(g_2/4) z +(g_3/4), \lb{1.3.78}  \\
&  E_{0}=-\wp(\Omega_1), \; E_{1}=-\wp(\Omega_2),\;
E_{2}=-\wp(\Omega_3). \no
\end{align}
Moreover, one has
\begin{align}
F_{1}(z,x) &= z+\wp(x+\Omega_3), \quad \mu_1(x)=-\wp(x+\Omega_3),
\lb{1.3.78a} \\ 
H_{2}(z,x) &= z^2-\wp(x+\Omega_3) z + \wp(x+\Omega_3)^2
-(g_2/4), \lb{1.3.79} \\
\nu_\ell(x)&=\big[\wp(x+\Omega_3)-(-1)^\ell[g_2
-3\wp(x+\Omega_3)^2]^{1/2}\big]{\big/}2, \quad \ell=0,1 \no 
\end{align}
and
\begin{align}
& \shKdV_1(V)=0, \lb{1.3.77b} \\
& \shKdV_2(V)-(g_2/8)\,\shKdV_0(V)=0, \text{ etc.} \lb{1.3.77c}
\end{align}
In addition, we record 
\begin{align}
&\psi_\pm(z,x,x_0)=\f{\sigma(x+\Omega_3\pm b)\sigma(x_0
+\Omega_3)}{\sigma(x+\Omega_3)\sigma(x_0+\Omega_3\pm b)} e^{\mp\zeta(b)(x-x_0)}, \\
&\psi_\pm(z,x+2\Omega_1,x_0)=\rho_\pm(z)\psi_\pm(z,x,x_0), \quad 
\rho_\pm(z)=e^{\pm[(b/\Omega_1)\zeta(\Omega_1)-\zeta(b)]2\Omega_1}
\end{align}
with Floquet parameter corresponding to $\Omega_1$-direction given by
\begin{equation}
k_1(b)=i[\zeta(b)\Omega_1-\zeta(\Omega_1) b]/\Omega_1. 
\end{equation}
Here 
\begin{align}
\begin{split}
&P=(z,y)=(-\wp(b),-(i/2)\wp'(b))\in\Pi_+, \lb{3.132} \\
&P^*=(z,-y)=(-\wp(b),(i/2)\wp'(b))\in\Pi_-,
\end{split}
\end{align}
where $b$ varies in the fundamental period parallelogram spanned
by the vertices $0$, $2\Omega_1$, $2\Omega_2$, and $2\Omega_3$. One then
computes 
\begin{align}
& \Delta(z)=\cosh[2(\Omega_1\zeta(b)-b\zeta(\Omega_1))], \\
& \langle \mu_1 \rangle = \zeta(\Omega_1)/\Omega_1, \quad 
\langle V\rangle = -2\zeta(\Omega_1)/\Omega_1, \\
& g(z,x)=-\f{z+\wp(x+\Omega_3)}{\wp'(b)}, \\ 
& \f{d}{dz}\big\langle g(z,\cdot)^{-1}\big\rangle
=2\f{z-[\zeta(\Omega_1)/\Omega_1]}{\wp'(b)}=-2\langle
g(z,\cdot)\rangle, \\ 
& \big\langle g(z,\cdot)^{-1}\big\rangle = -2 [\zeta(b)
-(b/\Omega_1)\zeta(\Omega_1)],
\end{align}
where $(z,y)=(-\wp(b),-(i/2)\wp'(b))\in\Pi_+$. The spectrum of the
operator $H$ with potential $V(x)=2\wp(x+\Omega_3)$ is then determined
as follows
\begin{align}
\sigma(H)&= \{\lambda\in\bbC\,|\, \Delta(\lambda)\in [-1,1]\} \\
&= \big\{\lambda\in\bbC\,\big|\, \Re\big(\big\langle
g(\lambda,\cdot)^{-1}\big\rangle\big)=0\big\} \\
&= \{\lambda\in\bbC\,|\,
\Re[\Omega_1\zeta(b)-b\zeta(\Omega_1)]=0, \, \lambda=-\wp(b)\}.
\end{align}
\end{example}

Generically (cf.\ \cite{Tk94}), $\sigma(H)$ consists of one simple
analytic arc (connecting two of the three branch points $E_m$, $m=0,1,2$)
and one simple semi-infinite analytic arc (connecting the remaining of the
branch points and infinity). The semi-infinite arc $\sigma_\infty$
asymptotically approaches the half-line
$L_{\langle V\rangle}=\{z\in\bbC\,|\, z=-2\zeta(\Omega_1)/\Omega_1 +x,
\, x\geq 0\}$ in the following sense: asymptotically,
$\sigma_\infty$ can be parameterized by
\begin{equation}
\sigma_\infty=\big\{z\in\bbC \,\big|\,
z=R-2i\,[\Im(\zeta(\Omega_1))/\Omega_1] +\Oh\big(R^{-1/2}\big) 
\text{ as $R\uparrow\infty$}\big\}. 
\end{equation}

We note that a slight change in the setup of Example \ref{eC.1}
permits one to construct crossing spectral arcs as shown in \cite{GW95}.
One only needs to choose complex conjugate fundamental half-periods
$\hatt \Omega_1\notin\bbR$, $\hatt \Omega_3=\ol{\hatt \Omega_1}$ with real
period $\Omega=2\big(\hatt \Omega_1+\hatt \Omega_3\big)>0$ and consider
the potential $V(x)=2\wp\big(x+a\,\big|\hatt \Omega_1,\hatt
\Omega_3\big)$, $0<\Im(a) <2\big|\Im\big(\hatt \Omega_1\big)\big|$.

Finally, we briefly consider a change of homology basis and illustrate
Theorem \ref{tB.5}. Let $\Omega_1>0$ and $\Omega_3\in\bbC$,
$\Im(\Omega_3)>0$. We choose the homology basis $\{\ti a_1,\ti b_1\}$
such that $\ti b_1$ encircles $E_0$ and
$E_1$ counterclockwise on $\Pi_+$ and $\ti a_1$ starts near $E_1$,
intersects $\ti b_1$ on $\Pi_+$, surrounds $E_2$ clockwise and then
continues on $\Pi_-$ back to its initial point surrounding $E_1$ such
that \eqref{aa16a} holds. Then,
\begin{align}
& \omega_1=c_1(1)\,dz/y, \quad c_1(1)=(4i\Omega_1)^{-1}, \\
& \int_{\ti a_1} \omega_1 =1, \quad \int_{\ti b_1} \omega_1 = \tau,
\quad \tau=\Omega_3/\Omega_1, \\
& \wti\omega^{(2)}_{\Pinf,0}=-\f{(z-\lambda_1)dz}{2y}, \quad 
\lambda_1=\zeta(\Omega_1)/\Omega_1, \\
& \int_{\ti a_1} \wti\omega^{(2)}_{\Pinf,0}=0, \quad 
\f{1}{2\pi i}\int_{\ti b_1}
\wti\omega^{(2)}_{\Pinf,0}=-2c_1(1)=\wti U_{0,1}, \\  
& \wti U_{0,1}=\f{i}{2\Omega_1}\in i \bbR, \\
& \int_{Q_0}^{P} \wti\omega^{(2)}_{\Pinf,0} -\ti e_0^{(2)}(Q_0)
\underset{b\to 0}{=} \f{i}{b}+\Oh(b) \no \\
& \hspace*{3.25cm} \underset{\zeta\to 0}{=}-\zeta^{-1} +\Oh(\zeta), 
\quad \zeta=\sigma/z^{1/2}, \, \sigma\in\{1,-1\},  \\
& \ti e_0^{(2)}(Q_0)
=-i[\zeta(b_0)\Omega_1-\zeta(\Omega_1)b_0]/\Omega_1, \\ 
& i\bigg[\int_{Q_0}^{P} \wti\omega^{(2)}_{\Pinf,0}
-\ti e_0^{(2)}(Q_0)\bigg]=
[\zeta(\Omega_1)b-\zeta(b)\Omega_1]/\Omega_1,
\\  & P=(-\wp(b),-(i/2)\wp'(b)), \; Q_0=(-\wp(b_0),-(i/2)\wp'(b_0)).
\no   
\end{align}
The change of homology basis (cf.\ \eqref{aa22A}--\eqref{aa22G}) 
\begin{align}
&\begin{pmatrix} \ti a_1\\ \ti b_1\end{pmatrix} \mapsto 
\begin{pmatrix} a_1'\\b_1'\end{pmatrix} = 
\begin{pmatrix} A & B \\ C & D \end{pmatrix} 
\begin{pmatrix} \ti a_1\\ \ti b_1\end{pmatrix} = 
\begin{pmatrix} A \ti a_1+B \ti b_1 \\ C a_1+ D b_1\end{pmatrix} ,  \\
& A, B, C, D \in \bbZ, \quad AD-BC=1,  
\end{align}
then implies
\begin{align}
& \omega_1'=\f{\omega_1}{A+B\tau}, \\
& \tau'=\f{\Omega_3'}{\Omega_1'}=\f{C+D\tau}{A+B\tau}, \\ 
& \Omega_1'= A\Omega_1+B\Omega_3, \quad \Omega_3'=C\Omega_1+D\Omega_3, \\
& \omega^{(2)\,\prime}_{\Pinf,0}=-\f{(z-\lambda_1')dz}{2y}, \quad 
\lambda^\prime_1=\lambda_1-\f{\pi iB}{2\Omega_1\Omega_1^\prime}, \\
& \int_{a^\prime_1} \omega^{(2)\,\prime}_{\Pinf,0}=0, \quad 
\f{1}{2\pi i}\int_{b^\prime_1}
\omega^{(2)\,\prime}_{\Pinf,0}=-\f{2c_1(1)}{A+B\tau}=U_{0,1}^\prime, \\
& U_{0,1}^\prime=\f{\wti U_{0,1}}{A+B\tau}=\f{i}{2\Omega_1^\prime}. 
\end{align}
Moreover, one infers
\begin{align}
& \psi_\pm(z,x+2\Omega_1^\prime,x_0)=\rho_\pm(z)^\prime\psi_\pm(z,x,x_0),
\no \\
& \rho_\pm(z)^\prime
=e^{\pm[(b/\Omega_1^\prime)(A\zeta(\Omega_1)+B\zeta(\Omega_3))
-\zeta(b)]2\Omega_1^\prime}
\end{align}
with Floquet parameter $k_1(b)^\prime$ corresponding to
$\Omega_1^\prime$-direction given by
\begin{equation}
k_1(b)^\prime=i\bigg[\zeta(b)\Omega_1-\zeta(\Omega_1) b+ 
\f{\pi iB}{2\Omega_1^\prime}b\bigg]\bigg/\Omega_1. 
\end{equation}

\vspace*{2mm}
{\bf Acknowledgments.}
F.\ G.\ is particularly indebted to Vladimir A.\ Marchenko for renewing 
his interest in the spectral theoretic questions addressed in this paper
and for the discussions we shared on this topic in June of 2000 at the
Department of Mathematical Sciences of the Norwegian University of Science
and Technology in Trondheim, Norway. 

We thank Helge Holden and Kwang Shin for many discussions on topics
related to this paper and Kwang Shin for a critical reading of our
manuscript. Moreover, we are indebted to Norrie Everitt and Igor Verbitsky
for pointing out the origin of Lemma
\ref{l3.10} and to J\"urgen Voigt for pointing out references \cite{Ar94}
and \cite{Ku98} to us. 



\begin{thebibliography}{10}
%
\bi{AS72} M.\ Abramowitz and I.\ A.\ Stegun,
{\em Handbook of Mathematical Functions}, Dover, New York, 1972.
%
\bi{Ap80} P.\ E.\ Appell, \textit{Sur la transformation des \'equations
diff\'erentielles lin\'eaires}, Comptes Rendus \textbf{91}, 211--214 
(1880).
%
\bi{Ar94} W.\ Arendt, {\it Gaussian estimates and interpolation of the
spectrum in $L^p$}, Diff. Integral Eqs. {\bf 7}, 1153--1168 (1994).
%
\bi{AE97} W.\ Arendt and A.\ F.\ M.\ ter Elst, {\it Gaussian estimates for
second order elliptic operators with boundary conditions}, J. Operator Th.
{\bf 38}, 87--130 (1997).
%
\bi{Be86} A.\ F.\ Beardon, {\it A Primer on Riemann Surfaces}, London
Math. Soc. Lecture Notes, Vol.\ 78, Cambridge University Press,
Cambridge, 1986.
%
\bi{BB83} E.\ F.\ Beckenbach and R.\ Bellman, {\it Inequalities}, 4th
printing, Springer, Berlin, 1983.
%
\bibitem{BBEIM94} E.~D.~Belokolos, A.~I.~Bobenko, V.~Z.~Enol'skii, 
A.~R.~Its, and V.~B.~Matveev, {\em Algebro-Geometric Approach to
Nonlinear Integrable Equations}, Springer, Berlin, 1994.
%
\bi{Be54} A.\ S.\ Besicovitch, {\it Almost Periodic Functions}, Dover, New
York, 1954.
%
\bi{Bi86a} B.\ Birnir, {\it Complex Hill's equation and the complex
periodic Korteweg-de Vries equations}, Commun. Pure Appl. Math. 
{\bf 39}, 1--49 (1986).
%
\bi{Bi86b} B.\ Birnir, {\it Singularities of the complex
Korteweg-de Vries flows}, Commun. Pure Appl. Math. {\bf 39}, 283--305
(1986).
%
\bi{Bo47} H.\ Bohr, {\it Almost Periodic Functions}, Chelsea, New York,
1947.
%
\bi{CL90} R.\ Carmona and J.\ Lacroix, {\it Spectral Theory
of Random Schr\"odinger Operators}, Birkh\"auser, Boston, 1990. 
%
\bi{CE70} R.\ S.\ Chisholm and W.\ N.\ Everitt, {\it On bounded integral
operators in the space of integrable-square functions}, Proc. Roy. 
Soc. Edinburgh Sect. A {\bf 69}, 199--204 (1970/71).
%
\bi{CEL99} R.\ S.\ Chisholm, W.\ N.\ Everitt, L.\ L.\ Littlejohn, {\it An
integral operator inequality with applications}, J. Inequal. Appl. 
{\bf 3}, 245--266 (1999). 
%
\bi{Co89} C.\ Corduneanu, {\it Almost Periodic Functions}, 2nd ed.,
Chelsea, New York, 1989.
%
\bi{Du75} B.\ A.\ Dubrovin, {\it Periodic problems for the Korteweg-de 
Vries equation in the class of finite-gap potentials}, Funct. Anal. Appl.
{\bf 9}, 215--223 (1975).
%
\bi{DMN76} B.\ A.\ Dubrovin, V.\ B.\ Matveev, and S.\ P.\ Novikov,
{\it Non-linear equations of Korteweg-de Vries type, finite-zone 
linear operators, and Abelian varieties}, Russian Math. Surv. 
{\bf 31:1}, 59--146 (1976).
%
\bi{Ea67} M.\ S.\ P.\ Eastham, {\it Gaps in the essential spectrum
associated with singular differential operators}, Quart. J. Math. {\bf
18}, 155--168 (1967).
%
\bi{Ea73} M.\ S.\ P.\ Eastham, {\it The Spectral Theory of Periodic 
Differential Equations}, Scottish Academic Press, Edinburgh and London,
1973.
%
\bi{EE89} D. E. Edmunds and W. D. Evans, {\it Spectral Theory and
Differential Operators}, Clarendon Press, Oxford, 1989. 
%
\bi{FK92} H.\ M.\ Farkas and I.\ Kra,
{\it Riemann Surfaces}, 2nd ed., Springer, New York, 1992.
%
\bi{Fi74} A.\ M. Fink, {\it Almost Periodic Differential Equations},
Lecture Notes in Math. {\bf 377}, Springer, Berlin, 1974. 
%
\bi{Fl75} H.\ Flaschka, {\it On the inverse problem for Hill's operator}, 
Arch. Rat. Mech. Anal. {\bf 59}, 293--309 (1975).
%
\bi{GD75}  I.\ M.\ Gel'fand and L.\ A.\ Dikii, {\it Asymptotic
behaviour of the resolvent of Sturm-Liouville equations and the
algebra of the Korteweg-de Vries equations}, Russ. Math. 
Surv. {\bf 30:5}, 77--113 (1975).
%
\bi{Ge01} F.\ Gesztesy, {\it Integrable systems in the infinite genus
limit}, Diff.\ Integral Eqs.\ {\bf 14}, 671--700 (2001).
%
\bi{GH03} F.\ Gesztesy and H.\ Holden, {\it Soliton Equations and
Their Algebro-Geometric Solutions. Vol. I: $(1+1)$-Dimensional 
Continuous Models},  Cambridge Studies in Advanced Mathematics,
Vol.\ 79, Cambridge Univ. Press, 2003.
%
\bi{GRT96} F.\ Gesztesy, R.\ Ratnaseelan, and G.\ Teschl.
{\em The KdV hierarchy and associated trace formulas}, in {\em Recent
 Developments in Operator Theory and Its Applications}, 
I.\ Gohberg, P.\ Lancaster, and P.\ N.\ Shivakumar, editors, 
{\em Operator Theory: Advances and Applications},  Vol.\ 87,  
Birkh\"auser, Basel, 1996, pp. 125--163.
%
\bi{GW95} F.\ Gesztesy and R.\ Weikard, {\it Floquet theory
revisited}, in {\it Differential Equations and Mathematical
Physics}, I.~Knowles (ed.), International Press, Boston, 1995, pp.\
67--84.
%
\bi{GW96} F.\ Gesztesy and R.\ Weikard, {\it Picard potentials and Hill's
equation on a torus}, Acta Math. {\bf 176}, 73--107 (1996).
%
\bi{GW98} F.\ Gesztesy and R.\ Weikard, {\it Elliptic
algebro-geometric solutions of the KdV and AKNS hierarchies -- an
analytic approach}, Bull. Amer. Math. Soc. {\bf 35}, 271--317 (1998).
%
\bi{Gl65} I.\ M. Glazman, {\it Direct Methods of Qualitative Spectral
Analysis of Singular Differential Operators}, Moscow, 1963. English
Translation by Israel Program for Scientific Translations, 1965.
%
\bi{Go85} S.\ Goldberg, {\it Unbounded Linear Operators}, Dover, New York,
1985.
%
\bi{IM75} A.\ R.\ Its and V.\ B.\ Matveev, {\it Schr\"odinger 
operators with finite-gap spectrum and $N$-soliton solutions of 
the Korteweg-de Vries equation}, Theoret. Math. Phys. {\bf 23}, 343--355
(1975).
%
\bi{JM82} R.\ Johnson and J.\ Moser, {\it The rotation number for
almost periodic potentials}, Commun. Math. Phys. {\bf 84}, 403--438 
(1982).
%
\bi{Ka80} T.\ Kato, {\it Perturbation Theory for Linear Operators},
corr. printing of the 2nd ed., Springer, Berlin, 1980.
%
\bi{Ko84} S.\ Kotani, {\it Ljapunov indices determine absolutely
continuous spectra of stationary random one-dimensional Schr\"odinger
operators}, in {\it Stochastic Analysis}, K.\ It{\v o} (ed.),
North-Holland, Amsterdam, 1984, p.\ 225--247. 
%
\bi{Ko85} S.\ Kotani, {\it On an inverse problem for random Schr\"odinger
operators}, Contemporary Math. {\bf 41}, 267--281 (1985).
%
\bi{Ko87a} S.\ Kotani, {\it One-dimensional random Schr\"odinger 
operators and Herglotz functions}, in {\it Probabilistic 
Methods in Mathematical Physics}, K.\ It{\v o} and N.\ Ikeda (eds.),
Academic  Press, New York, 1987, p.\ 219--250. 
%
\bi{Ko97} S.\ Kotani, {\it Generalized Floquet theory for stationary
Schr\"odinger operators in one dimension}, Chaos, Solitons \& Fractals 
{\bf 8}, 1817--1854 (1997).
%
\bi{Ku98} P.\ C.\ Kunstmann, {\it Heat kernel estimates and $L^p$ spectral
independence of elliptic operators}, Bull. London Math. Soc. {\bf 31},
345--353 (1998).
%
\bi{La75} P.\ D.\ Lax, {\it Periodic solutions of the Korteweg--de Vries
equation}, Commun. Pure Appl. Math. {\bf 28}, 141--188.
%
\bi{Le87} B.\ M.\ Levitan, {\it Inverse Sturm--Liouville 
Problems,} VNU Science Press, Utrecht, 1987.
%
\bi{LZ82} B.\ M.\ Levitan and V.\ V.\ Zhikov, {\it Almost Periodic
Functions and Differential Equations}, Cambridge University Press,
Cambridge, 1982.
%
\bi{Ma74} V.\ A.\ Marchenko, {\it A periodic Korteweg--de Vries problem},
Sov. Math. Dokl. {\bf 15}, 1052--1056 (1974).
%
\bi{Ma74a} V.\ A.\ Marchenko, {\it The periodic Korteweg--de Vries
problem}, Math. USSR Sbornik {\bf 24}, 319--344 (1974).
%
\bibitem{Ma86} V.\ A.\ Marchenko, {\it Sturm-Liouville
Operators and Applications}, Birkh\"auser, Basel, 1986.
%
\bi{Ma92} A.\ I.\ Markushevich, {\it Introduction to the Classical
Theory of Abelian Functions}, Amer. Math. Soc, Providence, 1992.
%
\bi{MM75}  H.\ P.\ McKean and P.\ van Moerbeke, 
{\it The spectrum of Hill's equation}, Invent. Math. {\bf 30}, 217--274  
(1975).
%
\bi{Mu72} B. Muckenhoupt, {\it Hardy's inequality with weights}, 
Studia Math. {\bf 44}, 31--38 (1972).
%
\bi{Mu84} D.\ Mumford, {\em Tata Lectures on Theta {II}}, Birkh\" auser,
Boston, 1984.
%
\bi{No74} S.\ P.\ Novikov, {\it The periodic problem for the
Korteweg-de Vries equation}, Funct. Anal. Appl. {\bf 8}, 236--246
(1974).
%
\bibitem{NMPZ84} S.\ Novikov, S.\ V.\ Manakov, L.\ P.\ Pitaevskii, and
V.\ E.\ Zakharov, {\it Theory of Solitons}, Consultants Bureau, New
York, 1984.
%
\bi{PT91} L.\ A.\ Pastur and V.\ A.\ Tkachenko, {\it Geometry of the
spectrum of the one-dimensional Schr\"odinger equation with a periodic
complex-valued potential}, Math. Notes {\bf 50}, 1045-1050 (1991).
%
\bibitem{RS78} M.\ Reed and B.\ Simon, {\it Methods of Modern Mathematical
Physics. IV: Analysis of Operators,} Academic Press, New York, 1978.
%
\bi{Ro63} F.\ S.\ Rofe-Beketov, {\it The spectrum of non-selfadjoint
differential operators with periodic coefficients}, Sov. Math. Dokl. {\bf
4}, 1563--1566 (1963).
%
\bi{Sc65} G.\ Scharf, {\it Fastperiodische Potentiale}, Helv. Phys. Acta,
{\bf 38}, 573--605 (1965).
%
\bi{Se60} M.\ I.\ Serov, {\it Certain properties of the spectrum of a
non-selfadjoint differential operator of the second kind}, Sov. Math.
Dokl. {\bf 1}, 190--192 (1960).
%
\bi{Si82} B.\ Simon, {\it Almost periodic Schr\"odinger operators: A
review}, Adv. Appl. Math. {\bf 3}, 463--490 (1982).
%
\bi{Ta69} G.\ Talenti, {\it Osservazioni sopra una classe di
disuguaglianze}, Rend. Sem. Mat. Fis. Milano {\bf 39}, 171--185 (1969).
%
%
\bi{Tk64} V.\ A.\ Tkachenko, {\it Spectral analysis of the
one-dimensional Schr\"odinger operator with a periodic complex-valued
potential}, Sov. Math. Dokl. {\bf 5}, 413--415 (1964).
%
\bi{Tk94} V.\ A.\ Tkachenko, {\it Discriminants and generic spectra of
non-selfadjoint Hill's operators}, Adv. Soviet Math. {\bf 19}, 41--71
(1994).
%
\bi{Tk96} V.\ Tkachenko, {\it Spectra of non-selfadjoint Hill's
operators and a class of Riemann surfaces}, Ann. of Math. {\bf 143},
181--231 (1996). 
%
\bi{To69} G.\ Tomaselli, {\it A class of inequalities}, Boll. Un. Mat.
Ital. {\bf 21}, 622--631 (1969). 
%
\bi{We98} R.\ Weikard, {\it Picard operators}, Math. Nachr. {\bf 195}, 
251--266 (1998).
%
\bi{We98a} R.\ Weikard, {\it On Hill's equation with a singular
complex-valued potential}, Proc. London Math. Soc. (3) {\bf 76}, 603--633
(1998).
%
\end{thebibliography}
\end{document}